\documentclass[reqno,12pt,a4paper]{amsart}
\usepackage[foot]{amsaddr}

\usepackage{mathrsfs}
\usepackage{amssymb}
\usepackage{amsfonts}
\usepackage{latexsym}
\usepackage{amsthm}

\usepackage{graphicx}
\def\lb{\label}

\newcommand{\er}[1]{\textrm{(\ref{#1})}}

\begin{document}


\renewcommand{\theequation}{\arabic{section}.\arabic{equation}}
\theoremstyle{plain}
\newtheorem{theorem}{\bf Theorem}[section]
\newtheorem{lemma}[theorem]{\bf Lemma}
\newtheorem{corollary}[theorem]{\bf Corollary}
\newtheorem{proposition}[theorem]{\bf Proposition}
\newtheorem{definition}[theorem]{\bf Definition}
\newtheorem{remark}[theorem]{\it Remark}

\def\a{\alpha}  \def\cA{{\mathcal A}}     \def\bA{{\bf A}}  \def\mA{{\mathscr A}}
\def\b{\beta}   \def\cB{{\mathcal B}}     \def\bB{{\bf B}}  \def\mB{{\mathscr B}}
\def\g{\gamma}  \def\cC{{\mathcal C}}     \def\bC{{\bf C}}  \def\mC{{\mathscr C}}
\def\G{\Gamma}  \def\cD{{\mathcal D}}     \def\bD{{\bf D}}  \def\mD{{\mathscr D}}
\def\d{\delta}  \def\cE{{\mathcal E}}     \def\bE{{\bf E}}  \def\mE{{\mathscr E}}
\def\D{\Delta}  \def\cF{{\mathcal F}}     \def\bF{{\bf F}}  \def\mF{{\mathscr F}}
\def\c{\chi}    \def\cG{{\mathcal G}}     \def\bG{{\bf G}}  \def\mG{{\mathscr G}}
\def\z{\zeta}   \def\cH{{\mathcal H}}     \def\bH{{\bf H}}  \def\mH{{\mathscr H}}
\def\e{\eta}    \def\cI{{\mathcal I}}     \def\bI{{\bf I}}  \def\mI{{\mathscr I}}
\def\p{\psi}    \def\cJ{{\mathcal J}}     \def\bJ{{\bf J}}  \def\mJ{{\mathscr J}}
\def\vT{\Theta} \def\cK{{\mathcal K}}     \def\bK{{\bf K}}  \def\mK{{\mathscr K}}
\def\k{\kappa}  \def\cL{{\mathcal L}}     \def\bL{{\bf L}}  \def\mL{{\mathscr L}}
\def\l{\lambda} \def\cM{{\mathcal M}}     \def\bM{{\bf M}}  \def\mM{{\mathscr M}}
\def\L{\Lambda} \def\cN{{\mathcal N}}     \def\bN{{\bf N}}  \def\mN{{\mathscr N}}
\def\m{\mu}     \def\cO{{\mathcal O}}     \def\bO{{\bf O}}  \def\mO{{\mathscr O}}
\def\n{\nu}     \def\cP{{\mathcal P}}     \def\bP{{\bf P}}  \def\mP{{\mathscr P}}
\def\r{\varrho} \def\cQ{{\mathcal Q}}     \def\bQ{{\bf Q}}  \def\mQ{{\mathscr Q}}
\def\s{\sigma}  \def\cR{{\mathcal R}}     \def\bR{{\bf R}}  \def\mR{{\mathscr R}}
\def\S{\Sigma}  \def\cS{{\mathcal S}}     \def\bS{{\bf S}}  \def\mS{{\mathscr S}}
\def\t{\tau}    \def\cT{{\mathcal T}}     \def\bT{{\bf T}}  \def\mT{{\mathscr T}}
\def\f{\phi}    \def\cU{{\mathcal U}}     \def\bU{{\bf U}}  \def\mU{{\mathscr U}}
\def\F{\Phi}    \def\cV{{\mathcal V}}     \def\bV{{\bf V}}  \def\mV{{\mathscr V}}
\def\P{\Psi}    \def\cW{{\mathcal W}}     \def\bW{{\bf W}}  \def\mW{{\mathscr W}}
\def\o{\omega}  \def\cX{{\mathcal X}}     \def\bX{{\bf X}}  \def\mX{{\mathscr X}}
\def\x{\xi}     \def\cY{{\mathcal Y}}     \def\bY{{\bf Y}}  \def\mY{{\mathscr Y}}
\def\X{\Xi}     \def\cZ{{\mathcal Z}}     \def\bZ{{\bf Z}}  \def\mZ{{\mathscr Z}}
\def\O{\Omega}

\newcommand{\gA}{\mathfrak{A}}
\newcommand{\gB}{\mathfrak{B}}
\newcommand{\gC}{\mathfrak{C}}
\newcommand{\gD}{\mathfrak{D}}
\newcommand{\gE}{\mathfrak{E}}
\newcommand{\gF}{\mathfrak{F}}
\newcommand{\gG}{\mathfrak{G}}
\newcommand{\gH}{\mathfrak{H}}
\newcommand{\gI}{\mathfrak{I}}
\newcommand{\gJ}{\mathfrak{J}}
\newcommand{\gK}{\mathfrak{K}}
\newcommand{\gL}{\mathfrak{L}}
\newcommand{\gM}{\mathfrak{M}}
\newcommand{\gN}{\mathfrak{N}}
\newcommand{\gO}{\mathfrak{O}}
\newcommand{\gP}{\mathfrak{P}}
\newcommand{\gR}{\mathfrak{R}}
\newcommand{\gS}{\mathfrak{S}}
\newcommand{\gT}{\mathfrak{T}}
\newcommand{\gU}{\mathfrak{U}}
\newcommand{\gV}{\mathfrak{V}}
\newcommand{\gW}{\mathfrak{W}}
\newcommand{\gX}{\mathfrak{X}}
\newcommand{\gY}{\mathfrak{Y}}
\newcommand{\gZ}{\mathfrak{Z}}

\def\ve{\varepsilon}   \def\vt{\vartheta}    \def\vp{\varphi}
\def\vk{\varkappa}

\def\Z{{\mathbb Z}}    \def\R{{\mathbb R}}   \def\C{{\mathbb C}}
\def\T{{\mathbb T}}    \def\N{{\mathbb N}}   \def\dD{{\mathbb D}}


\def\la{\leftarrow}              \def\ra{\rightarrow}            \def\Ra{\Rightarrow}
\def\ua{\uparrow}                \def\da{\downarrow}
\def\lra{\leftrightarrow}        \def\Lra{\Leftrightarrow}


\def\lt{\biggl}                  \def\rt{\biggr}
\def\ol{\overline}               \def\wt{\widetilde}
\def\no{\noindent}


\let\ge\geqslant                 \let\le\leqslant
\def\lan{\langle}                \def\ran{\rangle}
\def\/{\over}                    \def\iy{\infty}
\def\sm{\setminus}               \def\es{\emptyset}
\def\ss{\subset}                 \def\ts{\times}
\def\pa{\partial}                \def\os{\oplus}
\def\om{\ominus}                 \def\ev{\equiv}
\def\iint{\int\!\!\!\int}        \def\iintt{\mathop{\int\!\!\int\!\!\dots\!\!\int}\limits}
\def\el2{\ell^{\,2}}             \def\1{1\!\!1}
\def\sh{\sharp}
\def\wh{\widehat}

\def\all{\mathop{\mathrm{all}}\nolimits}
\def\where{\mathop{\mathrm{where}}\nolimits}
\def\as{\mathop{\mathrm{as}}\nolimits}
\def\Area{\mathop{\mathrm{Area}}\nolimits}
\def\arg{\mathop{\mathrm{arg}}\nolimits}
\def\const{\mathop{\mathrm{const}}\nolimits}
\def\det{\mathop{\mathrm{det}}\nolimits}
\def\diag{\mathop{\mathrm{diag}}\nolimits}
\def\diam{\mathop{\mathrm{diam}}\nolimits}
\def\dim{\mathop{\mathrm{dim}}\nolimits}
\def\dist{\mathop{\mathrm{dist}}\nolimits}
\def\Im{\mathop{\mathrm{Im}}\nolimits}
\def\Iso{\mathop{\mathrm{Iso}}\nolimits}
\def\Ker{\mathop{\mathrm{Ker}}\nolimits}
\def\Lip{\mathop{\mathrm{Lip}}\nolimits}
\def\rank{\mathop{\mathrm{rank}}\limits}
\def\Ran{\mathop{\mathrm{Ran}}\nolimits}
\def\Re{\mathop{\mathrm{Re}}\nolimits}
\def\Res{\mathop{\mathrm{Res}}\nolimits}
\def\res{\mathop{\mathrm{res}}\limits}
\def\sign{\mathop{\mathrm{sign}}\nolimits}
\def\span{\mathop{\mathrm{span}}\nolimits}
\def\supp{\mathop{\mathrm{supp}}\nolimits}
\def\Tr{\mathop{\mathrm{Tr}}\nolimits}
\def\BBox{\hspace{1mm}\vrule height6pt width5.5pt depth0pt \hspace{6pt}}


\newcommand\nh[2]{\widehat{#1}\vphantom{#1}^{(#2)}}
\def\dia{\diamond}

\def\Oplus{\bigoplus\nolimits}



\def\qqq{\qquad}
\def\qq{\quad}
\let\ge\geqslant
\let\le\leqslant
\let\geq\geqslant
\let\leq\leqslant
\newcommand{\ca}{\begin{cases}}
\newcommand{\ac}{\end{cases}}
\newcommand{\ma}{\begin{pmatrix}}
\newcommand{\am}{\end{pmatrix}}
\renewcommand{\[}{\begin{equation}}
\renewcommand{\]}{\end{equation}}
\def\bu{\bullet}

\title[{Resonances and eigenvalues for a 1D half-crystal}]
{On the resonances and eigenvalues for a 1D half-crystal with localised
impurity}

\date{\today}
\author[Evgeny L. Korotyaev]{Evgeny L. Korotyaev}
\address{Saint-Petersburg University;
 \ korotyaev@gmail.com}

\date{\today}
\author[Karl Michael Schmidt]{Karl Michael Schmidt}
\address{School of Mathematics, Cardiff University.
Senghennydd Road, CF24 4AG Cardiff, Wales UK;
 \ SchmidtKM@cf.ac.uk}



\begin{abstract}
\no We consider the Schr\"odinger operator $H$ on the half-line with a periodic
potential $p$ plus a compactly supported  potential $q$. For generic $p$, its
essential spectrum has an infinite sequence of open gaps.
We determine the asymptotics of the resonance counting function and show that,
for sufficiently high energy, each non-degenerate gap contains exactly one
eigenvalue or antibound state, giving asymptotics for their positions.
Conversely, for any potential $q$ and for any sequences
$(\s_n)_{1}^\iy, \s_n\in \{0,1\}$, and $(\vk_n)_1^\iy\in \ell^2, \vk_n\ge 0$,
there exists a potential $p$ such that $\vk_n$ is the length of the $n$-th gap,
$n\in\N$,
and $H$ has exactly $\s_n$ eigenvalues and $1-\s_n$ antibound state in each
high-energy gap.
Moreover, we show that between any two eigenvalues in a gap, there is an odd
number of antibound states, and hence deduce an asymptotic lower bound on the
number of antibound states in an adiabatic limit.

\end{abstract}

\maketitle

\section {Introduction and main results}
\setcounter{equation}{0}

Consider the Schr\"odinger operator $H$ acting in the Hilbert space $L^2(\R_+ )$ and  given by
$$
H=H_0+q,\qqq \qqq  H_0f=-f''+pf
$$
with the boundary condition $f(0)=0$.
Here $p$ is 1-periodic and $q$ is compactly supported and they satisfy
\[
\lb{0}
p\in L_{real}^1(\R/\Z),\qqq \qqq  q\in \cQ_t,
\]
where
$\cQ_t=\{q\in L_{real}^2(\R_+) \mid  \sup ( \supp (q))=t\}$
and we keep $t>0$ fixed throughout.
For later use, we set $n_t=\inf_{n\in \N, n\ge t} n$.
The spectrum of $H_0$ consists of an absolutely continuous part $\s_{ac}(H_0)=
\bigcup\limits_{n\in\N} \gS_n$ plus  at most one eigenvalue in each non-empty gap
$\g_n$, $n\in\N$, \cite{E}, \cite{Zh3}),
 where the bands $\gS_n$ and gaps $\g_n$ are given by (see Fig. 1)
$$
\gS_n=[E^+_{n-1},E^-_n],\ \ \qq \g_{n}=(E^-_{n},E^+_n)\qq (n\in\N),\qq
$$
and the $E_n$ satisfy
\[\lb{00}
E_0^+< E^-_1 \le E^+_1\dots\le E^+_{n-1}< E^-_n \le E^+_{n}<\dots
\]
It is known that there are infinitely many non-degenerate gaps, i.e.
$E_n^- < E_n^+$, unless $p$ is arbitrarily often differentiable \cite{Ho},
and all gaps are non-degenerate generically \cite{MO}, \cite{Si}.
Without loss of generality, we may assume $E_0^+=0$. The sequence (\er{00})
is the spectrum of the equation
\[
\lb{1}
-y''+p(x)y=\l y
\]
with the condition of 2-periodicity, $y(x+2)=y(x)$ $(x\in \R)$.
We also set $\g_0=(-\iy,E_0^+)$.
If a gap degenerates, $\g_n=\es $ for some $n$, then the corresponding bands $\gS_{n} $ and $\gS_{n+1}$ touch.
This happens when $E_n^-=E_n^+$; this number is then a double eigenvalue of the 2-periodic problem \er{1}. The lowest eigenvalue $E_0^+=0$ is always simple and has a 1-periodic eigenfunction.
Generally, the eigenfunctions corresponding to eigenvalues $E_{2n}^{\pm}$ are 1-periodic, those for $E_{2n+1}^{\pm}$ are 1-anti-periodic in the sense that
$y(x+1)=-y(x)$ $(x\in\R)$.

Consider the operator $\mH y=-y''+(p+q)y$ on the real line, where $p$ is periodic and and $q$ is compactly supported.
The spectrum of the operator $\mH$ consists of an absolutely continuous part
$\s_{ac}(\mH)=\s_{ac}(H_0)$ plus a finite number of simple eigenvalues in each
non-empty gap $\g_n, n\in\N_0 := \N\cup\{0\}$ (\cite{Rb}, \cite{F1}), and has
at most two
eigenvalues in every open gap $\g_n$ with sufficiently large $n$ (\cite{Rb}).
If $q_0:=\int_\R q(x)dx\ne 0$, then $\mH$ has precisely one eigenvalue (\cite{Zh1}, \cite{F2}, \cite{GS}) and one antibound state \cite{K4} in each non-empty gap $\g_n$ with sufficiently large $n$.
If  $q_0=0$, then there are, roughly speaking, either two eigenvalues and no antibound states or no eigenvalues and two antibound states
in each non-empty gap $\g_n$ with sufficiently large $n$ (\cite{K4}).
Similarly, the spectrum of $H$ defined in \er{0} consists of an absolutely
continuous part $\s_{ac}(H)=\s_{ac}(H_0)$ plus a finite number of simple
eigenvalues in each non-empty gap $\g_n$, $n\in\N_0$. This follows from the corresponding property of $\mH$ by the Glazman decomposition principle.

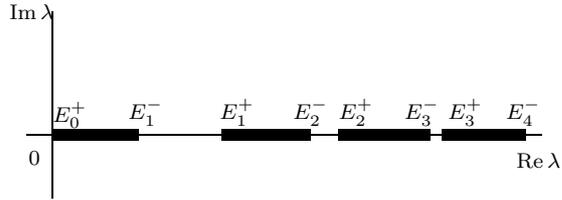
\begin{figure}
\tiny
\unitlength=1.00mm
\special{em:linewidth 0.4pt}
\linethickness{0.4pt}
\begin{picture}(108.67,33.67)
\put(41.00,17.33){\line(1,0){67.67}}
\put(44.33,9.00){\line(0,1){24.67}}
\put(108.33,14.00){\makebox(0,0)[cc]{$\Re\l$}}
\put(41.66,33.67){\makebox(0,0)[cc]{$\Im\l$}}
\put(42.00,14.33){\makebox(0,0)[cc]{$0$}}
\put(44.33,17.33){\linethickness{4.0pt}\line(1,0){11.33}}
\put(66.66,17.33){\linethickness{4.0pt}\line(1,0){11.67}}
\put(82.00,17.33){\linethickness{4.0pt}\line(1,0){12.00}}
\put(95.66,17.33){\linethickness{4.0pt}\line(1,0){11.00}}
\put(46.66,20.00){\makebox(0,0)[cc]{$E_0^+$}}
\put(56.66,20.33){\makebox(0,0)[cc]{$E_1^-$}}
\put(68.66,20.33){\makebox(0,0)[cc]{$E_1^+$}}
\put(78.33,20.33){\makebox(0,0)[cc]{$E_2^-$}}
\put(84.33,20.33){\makebox(0,0)[cc]{$E_2^+$}}
\put(93.00,20.33){\makebox(0,0)[cc]{$E_3^-$}}
\put(98.66,20.33){\makebox(0,0)[cc]{$E_3^+$}}
\put(106.33,20.33){\makebox(0,0)[cc]{$E_4^-$}}
\end{picture}
\caption{The cut domain $\C\sm \cup \gS_n$ and the cuts (bands) $\gS_n=[E^+_{n-1},E^-_n], n\ge 1$}
\lb{sS}
\end{figure}

Throughout the paper, we shall denote by
$\vt(x,z)$, $\vp(x,z)$ the two solutions forming the canonical fundamental
system of the unperturbed equation $-y''+py=z^2y$, i.e.,
satisfying the initial conditions $\vp'(0,z)=\vt(0,z)=1$ and $\vp(0,z)=\vt'(0,z)=0$. (Here and in the following $'$ denotes the derivative w.r.t. the first variable.)
The Lyapunov function (Hill discriminant) of the periodic equation is then
defined by $\D(z)={1\/2}(\vp'(1,z)+\vt(1,z))$.
The function $\l \mapsto \D^2(\sqrt \l)$ is entire, where we take the square
root to be positive on the positive real axis and to map the remainder of the
complex plane into the upper half-plane.
For the function $(1-\D^2(\sqrt \l))^{1\/2}$ $(\l\in \ol\C_+)$, we fix the
branch by the condition
$(1-\D^2(\sqrt {\l+i0}))^{1\/2}>0$ for $\l\in \gS_1=[E^+_{0},E^-_1]$, and
introduce the two-sheeted Riemann surface $\L$ of $(1-\D^2(\sqrt \l))^{1\/2}$
obtained by joining the upper and lower rims of two copies of the cut plane
$\C\sm\s_{ac}(H_0)$ in the usual (crosswise) way, see e.g. \cite{F3}.
We denote the $n$-th gap on the first, physical sheet $\L_1$ by $\g_n^{(1)}$
and its counterpart on the second, nonphysical sheet $\L_2$ by $\g_n^{(2)}$,
and set
\[
\lb{sL}
 \g_n^c:=\ol\g_n^{(1)}\cup \ol\g_n^{(2)}.
\]

It is well known (see e.g. \er{R})  that, for each
$h\in C_0^\iy(\R_+), h\ne 0$, the function
$f(\l)=((H-\l)^{-1}h,h)$ has a meromorphic extension from the physical sheet
$\L_1$ to the whole Riemann surface $\L$.
Moreover, if $f$ has a pole at some $\l_0\in \L_1$ for some $h$, then $\l_0$ is
an eigenvalue (bound state) of $H$ and
$\l_0\in\bigcup\limits_{n\in\N_0} \g_n^{(1)}$.

\bigskip
\no {\bf Definition. } Let $f(\l)=((H-\l)^{-1}h,h)$ $(\l\in \L)$ for some
$h\in C_0^\iy(\R_+), h\ne 0$.

\no 1) If  $f(\l)$ has a pole at some $\l_0\in \L_2$,   $\l_0\ne E_n^\pm, n\ge 0$, we call $\l_0$ a {\it resonance}.

\no 2)  A point  $\l_0=E_n^\pm , n\ge 0$, is called a {\it virtual state}
if the function $z \mapsto f(\l_0+z^2)$ has a pole at $0$.

\no 3) A point $\l_0\in\L$ is called a {\it state} if it is either a bound
state or a resonance or a virtual state.
Its {\it multiplicity} is the multiplicity of the corresponding pole.
We denote by $\gS_{st}(H)$ the set of all states.
If $\l_0\in \g_n^{(2)}, n\ge 0$, then we call $\l_0$ an {\it antibound state}.

\bigskip
\noindent
As a telling example we consider the states of the operator $H_0$ for the case
$p\ne \const , q=0$, see \cite{Zh3}, \cite{HKS}.
Let $f_0(\l)=((H_0-\l)^{-1}h,h)$ for some $h\in C_0^\iy(\R_+)$.
It is well known that the function  $f_0$ is meromorphic on the physical sheet
$\L_1$ and has a meromorphic extension into $\L$  (see e.g. \er{R}).
For each $\g_n^c\ne \es, n\ge 1$, there is exactly one state $\l_n^0\in \g_n^c$
of $H_0$ and its projection onto the complex plane coincides with the $n$-th
eigenvalue, $\m_n^2$, of  the Dirichlet
boundary value problem
$$
-y_n''+py_n=\m_n^2 y_n, \qq y_n(0)=y_n(1)=0, \qqq x\in [0, 1],\qq n\ge 1.
$$
Moreover, exactly one of the following three cases holds,

1) $\l_n^0\in \g_n^{(1)}$ is an eigenvalue,

2) $\l_n^0\in \g_n^{(2)}$ is an antibound state, or

3) $\l_n^0\in\{E_n^+, E_n^-\}$ is a virtual state.

\no There are no other states of $H_0$,
so $H_0$ has only eigenvalues, virtual states and antibound states.
If there are exactly $N\ge 1$ nondegenerate gaps in the spectrum of
$\s_{ac}(H_0)$, then the operator $H_0$ has exactly $N$ states; the closed gaps
$\g_n=\es$ do not contribute any states.
In particular, if $\g_n=\es$ $(n\ge 1)$, then $p=0$ (\cite{MO}, \cite{K5}) and $H_0$ has no states.
A more detailed description of the states of $H_0$ is given in Lemma \ref{Tm} below.

We shall need the following results from the inverse spectral theory for the
unperturbed operator $H_0$.
Defining the mapping $p\mapsto \x=(\x_n)_1^\iy$, where the components $\x_n=(\x_{1n},\x_{2n})\in \R^2$ are given by
$$
 \x_{1n}={E_n^-+E_n^+\/2}-\m_n^2,\qqq
\x_{2n}=\rt|{|\g_n|^2\/4}-\x_{1n}^2\rt|^{1\/2}a_n,\
\qq
a_n=\ca +1  & {\rm if} \ \l_n^0 \ {\rm is \ an \ eigenvalue,} \\
        -1 & {\rm if} \ \l_n^0 \ {\rm is \ a \ resonance,}      \\
         0 & {\rm if}  \ \l_n^0 \ {\rm is \ a \ virtual \ state,}  \ac
$$
we have the following result from \cite{K5}, \cite{K7}:

\no{\it The mapping $\x: \cH\to \ell^2\os \ell^2$ is a real analytic
isomorphism between the real Hilbert spaces
$\cH=\{p\in L^2(0,1)\mid \int_0^1p(x)dx=0\}$ and $\ell^2\os \ell^2$,
and the estimates
\[
\lb{esg}
\|p\|\le 4\|\x\| (1+ \|\x\|^{1\/3}),\qqq
\|\x\|\le \|p\|(1+\|p\|)^{1\/3}
\]
hold, where $\|p\|^2=\int_0^1p^2(x)dx$ and $\|\x\|^2={1\/4}\sum |\g_n|^2$.}

\no
Moreover, given any non-negative sequence $\vk=(\vk_n)_1^\iy\in \ell^2$,
there are unique 2-periodic eigenvalues $E_n^\pm$ $(n\in\N_0)$,
for some $p\in \cH$, such that each $\vk_n=E_n^+-E_n^-$ $(n\in\N)$.
Consequently, from the gap lengths $(|\g_n|)_1^\iy$ one can uniquely recover
the Riemann surface $\L$ as well as the points $E_n^-=E_n^+$ where $\vk_n=0$.
Furthermore, for any additional sequence $\wt\l_n^0\in \g_n^c$ $(n\in\N)$,
there is a unique potential $p\in \cH$ such that each state $\l_n^0$ of the
corresponding operator coincides with $\wt\l_n^0$ $(n\in\N)$.
The results of \cite{K5} were extended in \cite{K6} to periodic distributions
$p=w'$, where $w\in\cH$.

Now let
$$
D(\l)=\det (I+q(H_0-\l)^{-1}) \qqq (\l \in \C_+).
$$
The function
$D$ is analytic in $\C_+$ and has a meromorphic extension to $\L$.
Each zero of $D$ in $\L_1$ is an eigenvalue of $H$ and lies somewhere in the
union of the physical gaps $\bigcup\limits_{n\in\N_0} \g_n^{(1)}$.
So far only certain particular results are known concerning the zeros on the
non-physical sheet $\L_2$.
Note that the set of zeros of $D$ on $\L_2$ is symmetric with respect to the
real line, since $D$ is real on $\g_0^{(2)}$.

We now turn to the perturbed equation.
Let $\F (x,z)$ be the solution of the initial-value problem
\[
\lb{bcF}
-\F''+(p+q)\F=z^2 \F \ {\rm on}\ [0,\infty),\qqq \F (0,z)=0, \ \ \ \F'(0,z)=1
\qq (z\in\C).
\]
Then our first result is as follows.

\begin{theorem}
\lb{T1}

\no i)  Let $\Omega_\ve :=\{z\in\C\mid |z^2-E_n^\pm|\ge n\ve \ (n\in\N)\}$ for
some $\ve >0$.  Then the function $D$ satisfies
\[
\lb{T1-1}
D(z^2)=1+{\wh q(z)-\wh q(0)\/2iz}+{O(e^{t (|\Im z|-\Im z)})\/z^{2}}\
\ \ as \ \ |z|\to\iy, z\in\Omega_\ve 
\]
 where $\wh q(z)=\int_0^tq(x)e^{2izx}dx$, and
\[
\lb{st1}
\gS_{st}(H)\sm \gS_{st}(H_0)=\{ \l\in  \L\sm\gS_{st}(H_0)\mid\ D(\l)=0\}
\ss \L_2\cup \bigcup_{n\in\N_0}\g_n^{(1)},
\]
\[
\lb{st2}
\gS_{st}(H)\cap \gS_{st}(H_0)=\{ z\in \gS_{st}(H_0)\mid \F(n_t,z)=0 \}.
\]

\no ii) If $\l_n^0\in  \gS_{st}(H)\cap \gS_{st}(H_0)$, then
$\l_n^0\in \g_n^{c}\ne \es$ $(n\in\N)$ and
\[
\lb{ajf}
D(\l)\to D(\l_n^0)\ne 0\qq  as \qq \l\to \l_n^0.
\]
\no iii) (Logarithmic Law) Each resonance $\l\in \L_2$ of $H$  satisfies
\[
\lb{T1-2}
|\sqrt \l\sin \sqrt \l|\le C_Fe^{(2t+1)|\Im \sqrt \l|},\qqq
\ \ C_F=3(\|p\|_1+\|p+q\|_t)e^{2\|p+q\|_t+\|p\|_1},
\]
and there are no resonances in the domain
$\mD_{forb}=\{\l\in \L_2\sm \cup \ol\g_n^{(2)}\mid
4C_Fe^{2|\Im \sqrt \l|}<|\l|^{1\/2}\}$.
Here $\|f\|_s := \int_0^s |f(x)|\,dx$ $(s > 0)$.
\end{theorem}
\no {\bf Remarks.}
1) Let $\l_n^0\in \g_n^{(1)}$ be an eigenvalue of $H_0$ for some $n\ge 1$.
If $\F(n_t,\m_n)=0$, then $\m_n^2$ is an eigenvalue of the Dirichlet boundary
value problem $-y''+(p+q)y=\m_n^2 y, y(0)=y(n_t)=0$.
Then by \er{st2}, $\l_n^0$ is a bound state of $H$ and
\er{ajf} yields $D(\l_n^0)\ne 0$.
Thus $\l_n^0$ {\bf is a pole of a resolvent, but $\l_n^0$ is neither a zero of
$D$ nor a pole} of the S-matrix for $H,H_0$ given by
\[
\lb{sm}
\cS_M(z)={\ol {D(\l)\/ D(\l)}},\ \ \ \l \in \s_{ac}(H_0).
\]

\medskip
\no 2) If $D(\l)=0$ for some $\l=E_n^\pm\ne \m_n^2$, $n\in\N_0$, then by \er{st1}, $\l$ is a virtual state.

\medskip
\no 3) If $\m_n^2=E_n^\pm$ for some $n\in\N$, then by \er{st2},
$\m_n^2$ is a virtual state if and only if $\F(n_t,\m_n)=0$.

\medskip
\no 4) If $p=0$, then it is well known that each zero of $D$ is a state
(\cite{K1}, \cite{S}).
Moreover, each  resonance lies below a logarithmic curve depending only on $q$
(\cite{K1}, \cite{Z}).
The forbidden domain $\mD_{forb}\cap \C_-$ is similar to the one in the case
$p=0$, see \cite{K1}.

\begin{figure}
\tiny
\unitlength 1mm 
\linethickness{0.4pt}
\ifx\plotpoint\undefined\newsavebox{\plotpoint}\fi 
\begin{picture}(76.222,83.464)(0,0)
\put(-1.842,54.344){\line(1,0){78.064}}
\put(1.854,72.6){\line(0,-1){69.664}}
\qbezier(67.71,54.344)(69.782,56.136)(72.078,54.344)
\qbezier(53.71,54.344)(55.782,56.136)(58.078,54.344)
\qbezier(67.822,54.344)(69.39,53.056)(71.854,54.232)
\qbezier(53.822,54.344)(55.39,53.056)(57.854,54.232)
\qbezier(11.486,54.232)(16.974,56.696)(21.118,54.232)
\qbezier(11.486,54.232)(16.19,52.72)(20.894,54.12)
\qbezier(31.31,54.456)(34.95,56.696)(39.71,54.456)
\qbezier(30.974,54.232)(35.006,52.44)(39.486,54.232)
\bezier{50}(1.854,9.88)(21.622,10.888)(34.446,24.44)
\bezier{50}(34.446,24.328)(45.534,36.816)(45.87,54.456)
\qbezier(21.342,48.968)(35.398,24.776)(75.662,19.848)
\multiput(16.19,38.888)(.0331852,.0373333){27}{\line(0,1){.0373333}}
\multiput(19.662,43.144)(.0331852,.0373333){27}{\line(0,1){.0373333}}
\multiput(34.222,55.016)(.0331852,.0373333){27}{\line(0,1){.0373333}}
\multiput(54.606,54.344)(.0331852,.0373333){27}{\line(0,1){.0373333}}
\multiput(68.27,53.448)(.0331852,.0373333){27}{\line(0,1){.0373333}}
\multiput(37.134,54.904)(.0331852,.0373333){27}{\line(0,1){.0373333}}
\multiput(35.766,52.888)(.0331852,.0373333){27}{\line(0,1){.0373333}}
\multiput(26.046,46.056)(.0331852,.0373333){27}{\line(0,1){.0373333}}
\multiput(30.078,32.392)(.0331852,.0373333){27}{\line(0,1){.0373333}}
\multiput(44.19,24.216)(.0331852,.0373333){27}{\line(0,1){.0373333}}
\multiput(48.222,17.496)(.0331852,.0373333){27}{\line(0,1){.0373333}}
\multiput(57.518,21.416)(.0331852,.0373333){27}{\line(0,1){.0373333}}
\multiput(63.342,15.592)(.0331852,.0373333){27}{\line(0,1){.0373333}}
\multiput(71.07,18.952)(.0331852,.0373333){27}{\line(0,1){.0373333}}
\multiput(16.078,39.784)(.0373333,-.0331852){27}{\line(1,0){.0373333}}
\multiput(19.55,44.04)(.0373333,-.0331852){27}{\line(1,0){.0373333}}
\multiput(34.11,55.912)(.0373333,-.0331852){27}{\line(1,0){.0373333}}
\multiput(54.494,55.24)(.0373333,-.0331852){27}{\line(1,0){.0373333}}
\multiput(68.158,54.344)(.0373333,-.0331852){27}{\line(1,0){.0373333}}
\multiput(37.022,55.8)(.0373333,-.0331852){27}{\line(1,0){.0373333}}
\multiput(35.654,53.784)(.0373333,-.0331852){27}{\line(1,0){.0373333}}
\multiput(25.934,46.952)(.0373333,-.0331852){27}{\line(1,0){.0373333}}
\multiput(29.966,33.288)(.0373333,-.0331852){27}{\line(1,0){.0373333}}
\multiput(44.078,25.112)(.0373333,-.0331852){27}{\line(1,0){.0373333}}
\multiput(48.11,18.392)(.0373333,-.0331852){27}{\line(1,0){.0373333}}
\multiput(57.406,22.312)(.0373333,-.0331852){27}{\line(1,0){.0373333}}
\multiput(63.23,16.488)(.0373333,-.0331852){27}{\line(1,0){.0373333}}
\multiput(70.958,19.848)(.0373333,-.0331852){27}{\line(1,0){.0373333}}
\put(68.158,40.568){\makebox(0,0)[cc]{$the \ forbidden\ domain \
\mD_{forb}$}} \put(0.158,74.568){\makebox(0,0)[cc]{$ \Im \l$}}
\put(81.158,54.568){\makebox(0,0)[cc]{$ \Re \l$}}
\put(67.262,72.488){\makebox(0,0)[cc]{$the \ sheet \ \L_2 \ $}}
\put(3.15,51.72){\makebox(0,0)[cc]{$0$}}
\put(16.158,58.568){\makebox(0,0)[cc]{$ \g_1^1$}}
\put(16.158,51.568){\makebox(0,0)[cc]{$ \g_1^2$}}
\put(36.158,58.568){\makebox(0,0)[cc]{$ \g_2^1$}}
\put(36.158,51.568){\makebox(0,0)[cc]{$ \g_2^2$}}
\put(55.158,58.568){\makebox(0,0)[cc]{$ \g_3^1$}}
\put(55.158,51.568){\makebox(0,0)[cc]{$ \g_3^2$}}
\put(69.158,58.568){\makebox(0,0)[cc]{$ \g_4^1$}}
\put(69.158,51.568){\makebox(0,0)[cc]{$ \g_4^2$}}
\end{picture}
\caption{\footnotesize The eigenvalues on the physical gaps $\g_n^1\ss \L_1$, anti-bound states on the non-physical gaps $\g_n^2\ss \L_2$,
resonances on the non-physical sheet $\L_2$ and the forbidden domain $\mD_{forb}$} \lb{res}
\end{figure}
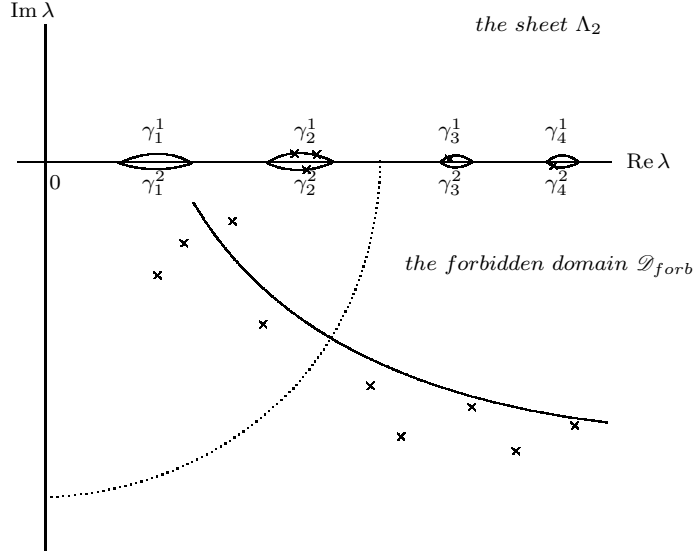

\bigskip
Let $\#(H,r, A)$ be the total number, counted according to multiplicity, of
states of $H$ of modulus $\le r$ in the set $A\subseteq \L$.

The Fourier coefficients $\wh p_{sn}, \wh q_{cn}$ and the Fourier transform
$\wh q$ are defined by
\[
\wh q_0=\!\int_0^t\!q(x)dx,\ \  \wh   p_{sn}=\!\int_0^1\!p(x)\sin 2\pi nxdx,\ \
\wh q(z)=\!\int_0^t\!q(x)e^{2izx}dx, \qq \wh q_{cn}=\Re \wh q(\pi n).
\]

\begin{theorem}
\lb{T2}
\no  i) $H$ has an odd number $\ge 1$ of states on each set
$\g_n^c\ne \es,n\ge 1$, where $\g_n^c$ is the union of the physical gap
$\ol \g_n^{(1)}\ss \L_1$ and the non-physical gap $\ol \g_n^{(2)}\ss\L_2$,
and $H$ has exactly one simple state $\l_n\in \g_n^c$ for all
$n\ge 1+4C_Fe^{t{\pi\/2}}$ with asymptotics
\[
\lb{T2-1}
\sqrt{\l_n}=\m_n-{(\wh q_{0}-\wh q_{cn})\wh p_{sn}\/2(\pi n)^2}+O({1\/n^3}) \qqq \as \qq n\to \iy.
\]
Moreover,  the following asymptotics hold true as $r\to \iy$:
\[
\lb{T2-2}
\#(H,r, \L_2\sm \cup \g_n^{(2)})=r{2t\/\pi}+o(r),
\]
\[
\lb{T2-3} \#(H,r, \g_\L)=\#(H_0,r, \g_\L)+2N_q \qqq for\ some\
integer \ N_q\ge 0, \qq r\notin \cup \ol\g_n,
\]
where $\g_\L=\bigcup_{n\ge 0} \ol\g_n^{(1)}\cup\ol\g_n^{(2)}$ is the
union of all gaps on $\L$.

 \no ii) Let $\l$ be an eigenvalue of $H$ and let $\l^{(2)}\in
\L_2$ be the same number but on the second sheet $\L_2$. Then
$\l^{(2)}$ is not an anti-bound state.

\no iii) Let $\l_1,\l_2\in \g_n^{(1)}, \l_1<\l_2,$ be eigenvalues of $H$  for some $n\ge 0$ and assume that there are no other eigenvalues on the interval $\O=(\l_1,\l_2)\ss \g_n^{(1)}$.
 Let $\O^{(2)}\ss \g_n^{(2)}\ss \L_2$ be the same interval but on the second sheet.
Then there exists an odd number $\ge 1$ of antibound states on  $\O^{(2)}$.

\end{theorem}

\no {\bf Remarks.}
1) Results (iii) with $p=0$ were obtained independently in \cite{K1}, \cite{S}.

\medskip
\no 2) The first term in the asymptotics \er{T2-2} is independent of the
periodic potential $p$. The asymptotics \er{T2-2} for the case $p=0$ was obtained by Zworski \cite{Z}.

\medskip
\no 3) The main difference between the distribution of the resonances in the
cases $p\ne \const$ and $p=\const$ concerns the bound states and antibound
states in high energy gaps, see \er{T2-1}.

\medskip
\no 4)  In the proof of \er{T2-2} we use a Paley-Wiener type theorem from
\cite{Fr}, the Levinson Theorem (see Sect. 4) and an analysis of the function
$D$ near $\l_n^0$.

\medskip
\no 5) For even potentials $p\in L_{even}^2(0,1)=\{p\in L^2(0,1)\mid
p(x)=p(1-x), x\in (0,1)\}$, all coefficients $p_{sn}$ vanish and the
asymptotics \er{T2-1} are not sharp.
This case is analyzed further in Theorem \ref{T4}.

\bigskip
We now turn to the question of stability of the real states $\l_n$.
As before, let $\l_n^0\in \g_n^c$ be a state of $H_0$.

\begin{theorem}
\lb{T3}
Let $b_n=q_0-\wh q_{cn}$ $(n\in\N)$.
Assume that $|p_{sn}|>n^{-\a}$ and $|b_n|>n^{-(1-\a)}$
for some $\a\in (0,1)$ and for all $n\in \cN_0$, where $\cN_0\ss \N$ is some
infinite subset such that $|\g_n|>0$ $(n\in \cN_0)$.
Assume $b_n>0$ (or $b_n<0$) for all $n\in\cN_0$.

Then each real state $\l_n\in\g_n^c$ of the perturbed operator $H$ has, for
sufficiently large $n\in \cN_0$, the following property.

\no If $\l_n^0$ is an eigenvalue of $H_0$, then $\l_n$ is an eigenvalue of $H$ and $\l_n^0<\l_n$ (or $\l_n^0>\l_n$).

\no If $\l_n^0$ is an antibound state of $H_0$, then $\l_n$ is an antibound
state of $H$ and $\l_n^0>\l_n$ (or $\l_n^0<\l_n$).

\end{theorem}

\no
{\bf Remarks.}
1) Let $q>0$. It is well known that the eigenvalues of $H_0+\t q$ are increasing
with the coupling constant $\t>0$. Roughly speaking, in the case considered in
Theorem \ref{T3} the antibound states in the gap move in the opposite direction.

\medskip
\no 2)
Numerical observations suggest the following scenario
as the coupling constant of the perturbation changes.
Consider the operator
$H_\t=H_0+\t q$, where $\t\in \R$ is  the coupling constant.
For $\t=0$, $H_0$ only has states $\l_n^0$, $n\in\N$ (eigenvalues, antibound
states and virtual states).
Take for instance the first gap $\g_1^c\ne \es$ and assume that $\l_1 = \l_1^0$
is an antibound state.
As $\t$ increases, the state $\l_1$ moves, and initially there are no other
states on $\g_1^c$. As $\t$ increases further, $\l_1$ eventually emerges onto
the physical gap $\g_1^{(1)}$ and becomes an eigenvalue; at first
there are no additional eigenvalues, but a pair of complex resonances
($\l\in \C_+\ss\L_2$ and $\ol \l\in \C_-\ss\L_2$) reach the non-physical gap
$\g_1^{(2)}$ and there turn into antibound states. With further increase of
$\t$, these antibound states will become virtual states and then bound states.
Thus, as $\t$ runs through $\R_+$, states not only move, but are also subject
to the following transmutations: resonances turn into antibound states, these
via virtual states into bound states, these via virtual states into antibound
states etc.

\bigskip
\noindent
We now consider the perturbations of virtual states.
Note that $p\in L_{even}^2(0,1)$  if and only if $\m_n^2 \in \{E_n^-, E_n^+\}$
 for all $n\in\N$ (\cite{GT}, \cite{KK1}).

\begin{theorem}
\lb{T4}
Let $p\in L^1(0,1)$ and assume there are unperturbed states
$\l_n^0\in \{E_n^-, E_n^+\}$ for all $n\in \cN_0$, where $\cN_0\ss \N$ is some
infinite subset such that $|\g_n|>0$ $(n\in \cN_0)$. Then
the following asymptotics for states of $H$ hold true:
\[
\lb{T4-1}
\sqrt{\l_n}=\m_n+s_n|\g_n|{(q_0-\wh q_{cn}+O({1\/n}))^2\/(2\pi n)^2}, \qqq
\qqq s_n=\ca +  & {\it if} \ \m_n^2=E_n^- \\
        - & {\it if} \ \m_n^2=E_n^+  \ac,\qq n\in \cN_0
\]
as   $n\to \iy$. Moreover, if  $\a\in ({1\/2},1)$, then for sufficiently large
$n\in\cN_0$, the following holds:

if $\l_n^0=E_n^-,\  q_0-\wh q_{cn}>n^{-\a}$ or $\l_n^0=E_n^+, \ q_0-\wh q_{cn}<-n^{-\a}$, then $\l_n$ is an eigenvalue,

if $\l_n^0=E_n^-, \ q_0-\wh q_{cn}<-n^{-\a}$ or  $\l_n^0=E_n^+, \
q_0-\wh q_{cn}>n^{-\a}$, then $\l_n$ is an antibound state.

\end{theorem}

{\bf Remark.} 1) Roughly speaking, \er{T4-1} gives the asymptotics for even
potentials $p\in L^1(0,1)$.

\medskip
\no 2) Consider the operator $-y''+(p+q+u)y$, $y(0)=0$, with
 an additional  potential perturbation
 $u\in L^2(\R_+)$ which is compactly
supported in $(0,t)$  and satisfies $|\wh u_{n}|=o(n^{-1})$ as
$n\to \iy$.
Then  the operator $H+u$ has the same number of bound states as $H$ in each gap
$\g_n\ne \es$ for $n$ large enough.

\bigskip

We have the following result on the inverse problem for our operator $H$.

\begin{theorem}
\lb{TIP}
i)  Let $q\in \cQ_t$ satisfy $|q_{0}-\wh q_{cn}|>n^{-\a}$ for sufficiently
large $n$, with some $\a\in ({1\/2},1)$.
Then for any sequences $(\s_n)_{1}^\iy, \s_n\in \{0,1\}$, and
$(\vk_n)_1^\iy\in\ell^2, \vk_n\ge 0$, there exists a potential $p\in L^2(0,1)$
such that the corresponding gap lengths $|\g_n|$ satisfy $|\g_n|=\vk_n, n\ge 1$,
and $H$ has exactly $\s_n$ eigenvalues and $1-\s_n$ antibound states in each gap
$\g_n\ne \es$ for sufficiently large $n$.

\no ii) Let $p\in L^1(0,1)$ and assume there are unperturbed states
$\l_n^0\in \{E_n^-, E_n^+\}$ for all $n\in \cN_0$, where $\cN_0\ss \N$ is some
infinite subset such that $|\g_n|>0$ $(n\in \cN_0)$. Then for any sequence $(\s_n)_{1}^\iy, \s_n\in \{0,1\}$, there exists a
potential $q\in \cQ_t$ such that $H$ has exactly $\s_n$ eigenvalues and
$1-\s_n$ antibound states in each gap $\g_n\ne \es$ for sufficiently large
$n\in \cN_0$.

\end{theorem}

\bigskip
\noindent
Let $\#_{bs}(H,\O)$ ($\#_{abs}(H,\O)$) be the total number, counted according to
multiplicity, of bound (antibound) states of $H$ on the segment
$\O\ss \ol\g_n^{(1)}\ss\L_1$ ($\O\ss \ol\g_n^{(2)}\ss\L_2$) for some $n\in\N_0$.

The integrated density of states $\r$ can be characterized as a
continuous, real-valued function on $\R$ with the properties
\[
\lb{ids}
\r([E_n^-,E_n^+])=n, \qqq \r(\gS_{n+1})=[n,n+1], \qqq
\cos \pi \r(\l)=\D(\sqrt \l) \ (\l\in \gS_{n+1})
\]
for all $n\in\N_0$.
The function $\r$ is strictly increasing on each spectral band $\gS_n$ and
constant on the closure of each gap $[E_n^-,E_n^+]$.
It is closely related to the quasimomentum defined in Section 2 via
$\r(\l)={1\/\pi}\Re k(\sqrt {\l+i0})$ $(\l\in \R)$.
Theorem \ref{T2} yields the following corollary.

\begin{theorem}
\lb{T5}
Let $H_\t=H_0+q_\t$ where $q_\t=q({\cdot\/\t}), \t\ge 1$.
Let $\O=[E_1,E_2]\ss \ol\g_n^{(1)}\ne \es$ be some interval on the physical
sheet $\L_1$ for some $n\ge 0$ and let $\O^{(2)}\ss \ol\g_n^{(2)}$ be the
corresponding interval on the non-physical sheet $\L_2$. Then
\begin{multline}
\lb{T5-1}
\#_{abs}(H_\t,\O^{(2)})\ge 1+\#_{bs}(H_\t,\O)
\\
=\t \int_0^\iy \rt(\r(E_2-q(x))-\r(E_1-q(x))\rt)dx+o(\t)\qq
\as \qqq \t\to \iy. \qq
\end{multline}
\end{theorem}

\no {\bf Remark.} 1) The proof of \er{T5-1} will be based on Sobolev's idea
\cite{So} for obtaining asymptotics of $\#_{bs}(H_\t,\O)$.
He considered the case $H_\t=H_0+\t V$, where $V(x)={c+o(1)\/x^\a}$ as
$x\to\iy$, for some $c\ne 0,\a>0$. This result is not directly applicable to
our case of a compactly supported perturbation. We therefore follow the
reinterpretation of the coupling constant as a scaling parameter, as previously
used in the study of perturbations of the periodic Dirac operator \cite{Sc},
which allows compactly supported perturbations.

\medskip
\no 2) Symmetry of bound states and antibound states for $p=0$
 in the semiclassical limit was studied in \cite{BZ}, \cite{DG}.

\bigskip
\noindent
A large number of papers are devoted to resonances for the
Schr\"odinger operator with $p=0$, see \cite{Fr}, \cite{H}, \cite{K1},
\cite{K2}, \cite{S}, \cite{Z}  and references therein.
Although resonances have been studied in many settings, there are relatively
few cases in which the asymptotics of the resonance counting function are
known, mainly in the one-dimensional setting (\cite{Fr}, \cite{K1}, \cite{K2},
\cite{S}, and \cite{Z}).
Zworski [Z] obtained the first results about the distribution of
resonances for the Schr\"odinger operator with compactly supported
potentials on the real line.
One of the present authors obtained the characterization, including uniqueness
and recovery, of the $S$-matrix for the Schr\"odinger operator with a compactly
supported potential on the real line \cite{K2} and on the half-line \cite{K1},
see also \cite{Z1}, \cite{BKW} concerning the uniqueness.

The following stability results for the Schr\"odinger operator on the half line
can be found in \cite{K3}; here the real Hilbert spaces $\ell_\ve^2$ are given
by
$$
 \ell^2_{\ve}=\rt\{f=(f_n)_1^{\iy }\mid \  f_n\in \R,\qq
 \| f \|_{\ve}^2=\sum _{n\ge 1 }(2\pi n)^{2\ve}f_n^2 <\infty \rt\}
 \ \ (\ve\ge 0), \ \qqq \ell^2_0=\ell^2.
$$

\noindent
(i) If $\vk=(\vk)_1^\iy$ is a sequence of poles (eigenvalues
and resonances) of the S-matrix for some real compactly
supported potential $q$ and $\wt\vk-\vk\in\ell_\ve^2$ for some
$\ve>1$, then $\wt\vk$ is a sequence of poles
of the S-matrix for some uniquely determined real-valued compactly supported
potential $\wt q$.

\noindent
(ii) The measure associated with the poles of the S-matrix is the Carleson
measure, and the sum $\sum (1+|\vk_n|)^{-\a}, \a>1$, can be estimated in terms
of the $L^1$-norm of the potential $q$.

\medskip
Brown and Weikard \cite{BW} considered  the Schr\"odinger operator
$-y''+(p_A+q)y$ on the half-line, where $p_A$ is an algebro-geometric potential
and $q$ is a compactly supported perturbation.
They proved that the zeros of the Jost function determine $q$ uniquely.

Christiansen \cite{Ch} considered resonances associated with the Schr\"odinger
operator $-y''+(p_{S}+q)y$ on the real line, where $p_S$ is a step potential.
She determined the asymptotics of the resonance-counting function and showed
that the resonances determine $q$ uniquely.

The recent paper \cite{K4} gives the following results about the operator $\mH
=\mH_0+q$, $\mH_0=-{d^2\/dx^2}+p$ on the real line, where $p$ is periodic
and $q$ is compactly supported:
1) asymptotics of the resonance-counting function are determined,
2) a forbidden domain for the resonances is specified,
3) the asymptotics of eigenvalues and antibound states are determined,
4) for any sequence $(\s)_1^\iy, \s_n\in \{0,2\}$, there exists a compactly
supported potential $q$ such that $\mH$ has $\s_n$ bound states and $2-\s_n$
antibound states in each gap $\g_n\ne \es$ for $n$ large enough,
5) for any $q$ (with $q_0=0$) and for any sequences
$(\s_n)_{1}^\iy, \s_n\in \{0,2\}$, and   $(\vk_n)_1^\iy\in \ell^2, \vk_n\ge 0$,
there exists a potential $p\in L^2(0,1)$ such that gap $\g_n$ has length
$|\g_n|=\vk_n$ $(n\in\N)$ and $\mH$ has exactly $\s_n$ eigenvalues and $2-\s_n$ antibound
states in each gap $\g_n\ne \es$ for sufficiently large $n$.

Thus, comparing the results for $\mH$ on $\R$ and $H$ on $\R_+$, we find that
their properties are close for even potentials $p\in L_{even}^2(0,1)$,
since in this case the unperturbed operator $H_0$ has only virtual states;
however, if $p$ is not even, then the unperturbed operator $H_0$
has eigenvalues, virtual states and antibound states in general, but the
operator $\mH_0$ has exactly two virtual states in each non-empty gap.
This leads to the different properties of the perturbed operators $H, \mH$, and,
roughly speaking, the case of $H$ is more complicated since the set of states
of the unperturbed operator $H_0$ is already more complicated.

For related results in the multidimensional case, see \cite{G}, \cite{D} and
references therein.

\bigskip
The structure of the present paper is as follows.
In Section 2 we define the Riemann
surface associated with the momentum variable $z=\sqrt \l, \l\in \L$, and
describe preliminary results about fundamental solutions.
In Section 3 we study the states of $H$.
In Section 4 we prove the main Theorems \ref{T1}-\ref{T5}.
The proofs use properties of the quasimomentum,
a priori estimates from \cite{KK}, \cite{MO},   and
results from the inverse theory for the Hill operator \cite{K5}.
The analysis of the entire function
$F(z)=\vp(1,z)D(z)\ol D(z), z^2\in \s_{ac}(H)$,
plays an important role, since its zeros are, as we show, closely related to
the states.
Thus the spectral analysis of $H$ is reduced to a problem of entire
function theory.


\section {Preliminaries }
\setcounter{equation}{0}

In the following, we shall treat the momentum $z=\sqrt \l$ (as opposed to the
energy $\l$) as the principal spectral variable. Note that the functions
$\vp(x,z), \vt(x,z)$ are entire in $z\in \C$ (cf.\ \cite{Ea} Thm 1.7.2), and so
are the functions $\wt\vt,
\wt\vp$ defined in \er{eqfwt} below. The other functions considered in the
present paper are combinations of $\vp(x,z), \vt(x,z), \wt\vp(x,z), \wt\vt(x,z)$
and $m_\pm={\b\pm i\sin k\/ \vp(1,\cdot)}$. Thus we are led to consider the
Riemann surface $\cM$ corresponding to the analytic continuation of the
function $\sin k$ (cf.\ \cite{K4}).
Take the cut domain (see Fig.3)
\[
\lb{2}
\cZ=\C\sm \cup \ol g_n,\qqq {\rm where} \qq
g_n=(e_n^-,e_n^+)=-g_{-n},\qq e_n^\pm=\sqrt{E_n^\pm}>0,\qq
n\ge 1,
\]
and note that $\D(e_{n}^{\pm})=(-1)^n$.
If $\l\in \g_n,  n\ge 1$, then $z\in g_{\pm n}$, and if $\l\in \g_0=(-\iy,0)$,
then  $z\in g_0^\pm=i\R_\pm$. Slitting the $n$-th momentum gap  $g_n$
(if non-empty), we obtain a cut $g_n^c$ with an upper rim $g_n^+$ and lower rim
$g_n^-$. Below we will identify this cut $g_n^c$ and the union of the upper rim
 (gap) $\ol g_{n}^+$ and the lower rim (gap) $\ol g_{n}^{\ -}$, i.e.,
\[
g_n^c=\ol g_{n}^+\cup \ol g_{n}^-,\ {\rm where}\ g_{n}^\pm =g_n\pm i0;
\qq {\rm and}\ z\in g_n \Rightarrow z\pm i0\in g_n^\pm.
\]

\begin{figure}
\tiny
\unitlength=1mm
\special{em:linewidth 0.4pt}
\linethickness{0.4pt}
\begin{picture}(120.67,34.33)
\put(20.33,21.33){\line(1,0){100.33}}
\put(70.33,10.00){\line(0,1){24.33}}
\put(69.00,19.00){\makebox(0,0)[cc]{$0$}}
\put(120.33,19.00){\makebox(0,0)[cc]{$\Re z$}}
\put(67.00,33.67){\makebox(0,0)[cc]{$\Im z$}}
\put(81.33,21.33){\linethickness{2.0pt}\line(1,0){9.67}}
\put(100.33,21.33){\linethickness{2.0pt}\line(1,0){4.67}}
\put(116.67,21.33){\linethickness{2.0pt}\line(1,0){2.67}}
\put(60.00,21.33){\linethickness{2.0pt}\line(-1,0){9.33}}
\put(40.00,21.33){\linethickness{2.0pt}\line(-1,0){4.67}}
\put(24.33,21.33){\linethickness{2.0pt}\line(-1,0){2.33}}
\put(81.67,24.00){\makebox(0,0)[cc]{$e_1^-$}}
\put(91.00,24.00){\makebox(0,0)[cc]{$e_1^+$}}
\put(100.33,24.00){\makebox(0,0)[cc]{$e_2^-$}}
\put(105.00,24.00){\makebox(0,0)[cc]{$e_2^+$}}
\put(115.33,24.00){\makebox(0,0)[cc]{$e_3^-$}}
\put(120.00,24.00){\makebox(0,0)[cc]{$e_3^+$}}
\put(59.33,24.00){\makebox(0,0)[cc]{$-e_1^-$}}
\put(50.67,24.00){\makebox(0,0)[cc]{$-e_1^+$}}
\put(40.33,24.00){\makebox(0,0)[cc]{$-e_2^-$}}
\put(34.67,24.00){\makebox(0,0)[cc]{$-e_2^+$}}
\put(26.00,24.00){\makebox(0,0)[cc]{$-e_3^-$}}
\put(19.50,24.00){\makebox(0,0)[cc]{$-e_3^+$}}
\end{picture}
\caption{The cut domain $\cZ=\C\sm \cup \ol g_n$ and the slits $g_n=(e_n^-,e_n^+)$ in the $z$-plane.}
\lb{z}
\end{figure}
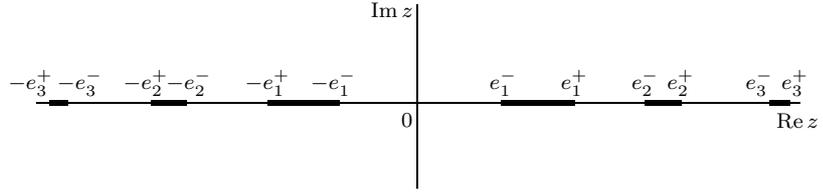

Now the function $\sin k(z)=\sqrt{1-\D^2(z)}, z\in \cZ$, which is analytic
because $\vp(x,z), \vt'(x,z)$ are, has an analytic
continuation through the cuts $\cup g_n\ss \R$, since $\sqrt{1-\D^2(z)}\in i\R$
for all $z\in \cup g_n$. Thus we obtain the analytic function $\sqrt{1-\D^2(z)}$
on the Riemann surface $\cM$, and (by applying the square root function) the
analytic function $\sqrt{1-\D^2(\sqrt \l)}$ on the Riemann surface $\L$.

The Riemann surface $\cM$ can be obtained from the cut domain
$\cZ=\C\sm \cup \ol g_n$ as follows. We identify the upper rim $g_{n}^+$ of the
cut $g_n^c$ with the upper rim $g_{-n}^+$ of the slit $g_{-n}^c$ and,
correspondingly,
the lower rim $g_{n}^-$ of the  cut $g_{n}^c$ with the lower rim $g_{-n}^-$ of
the  cut $g_{-n}^c$ for all non-empty gaps.
The mapping $z=\sqrt \l$ gives a bijection from $\L$ onto $\cM$.
The gap $\g_n^{(1)}\ss \L_1$ is mapped onto $g_n^+\ss \cM$ and the gap
$\g_n^{(2)}\ss \L_2$ is mapped onto $g_n^-\ss \cM$.
The upper rim $g_{n}^+$ corresponds to a physical gap, the lower rim $g_{n}^-$
to a non-physical gap.
Moreover, the union of $\cM\cap\C _+=\cZ\cap\C _+$ with all physical gaps
$g_{n}^+$ forms the so-called physical sheet $\cM_1$, while the union of
$\cM\cap\C _-=\cZ\cap\C _-$ with all non-physical gaps $g_{n}^-$ forms the
so-called non-physical sheet $\cM_2$. The two sheets are joined together on
the set $\R\sm \cup g_n$ (which is the spectrum of the periodic operator $H_0$).

\begin{figure}
\tiny
\unitlength=1mm
\special{em:linewidth 0.4pt}
\linethickness{0.4pt}
\begin{picture}(120.67,34.33)
\put(20.33,20.00){\line(1,0){102.33}}
\put(71.00,7.00){\line(0,1){27.00}}
\put(70.00,18.67){\makebox(0,0)[cc]{$0$}}
\put(124.00,18.00){\makebox(0,0)[cc]{$\Re k$}}
\put(67.00,33.67){\makebox(0,0)[cc]{$\Im k$}}
\put(87.00,15.00){\linethickness{2.0pt}\line(0,1){10.}}
\put(103.00,17.00){\linethickness{2.0pt}\line(0,1){6.}}
\put(119.00,18.00){\linethickness{2.0pt}\line(0,1){4.}}
\put(56.00,15.00){\linethickness{2.0pt}\line(0,1){10.}}
\put(39.00,17.00){\linethickness{2.0pt}\line(0,1){6.}}
\put(23.00,18.00){\linethickness{2.0pt}\line(0,1){4.}}
\put(85.50,18.50){\makebox(0,0)[cc]{$\pi$}}
\put(54.00,18.50){\makebox(0,0)[cc]{$-\pi$}}
\put(101.00,18.50){\makebox(0,0)[cc]{$2\pi$}}
\put(36.00,18.50){\makebox(0,0)[cc]{$-2\pi$}}
\put(117.00,18.50){\makebox(0,0)[cc]{$3\pi$}}
\put(20.00,18.50){\makebox(0,0)[cc]{$-3\pi$}}
\put(87.00,26.00){\makebox(0,0)[cc]{$\pi+ih_1$}}
\put(56.00,26.00){\makebox(0,0)[cc]{$-\pi+ih_1$}}
\put(103.00,24.00){\makebox(0,0)[cc]{$2\pi+ih_2$}}
\put(39.00,24.00){\makebox(0,0)[cc]{$-2\pi+ih_2$}}
\put(119.00,23.00){\makebox(0,0)[cc]{$3\pi+ih_3$}}
\put(23.00,23.00){\makebox(0,0)[cc]{$-3\pi+ih_3$}}
\end{picture}
\caption{The domain $\cK=\C\sm \cup \G_n$ with $\G_n=(\pi n-ih_n,\pi n+ih_n)$}
\lb{k}
\end{figure}
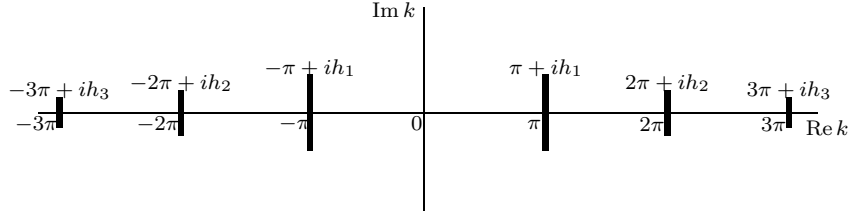

We introduce the {\it quasimomentum} $k$ for $H_0$ as $k(z)=\arccos \D(z)$
$(z \in \cZ)$. The function $k$ is analytic in $\cZ$ and satisfies
\[
\lb{pk}
(i)\qq k(z)=z+O(1/z)\qq  as \ \ |z|\to \iy, \qq  \qq (ii)\qq  \Re k(z\pm i0) |_{[e_n^-,e_n^+]}=\pi n\qq (n\in \Z),
\]
as well as $\pm \Im k(z)>0$ for any $z\in \C_\pm$ (\cite{MO},
\cite{KK}). The function sin of the  quasimomentum $k$ is also
analytic on $\cM$ and satisfies $\sin k(z)=(1-\D^2(z))^{1\/2}$ $(z\in \cM)$.
Moreover, $k$ is a conformal mapping from $\cZ$ onto the
quasimomentum domain $\cK=\C\sm \cup \G_n$, see  Figs. 3 and 4. Here
$\G_n=(\pi n-ih_n,\pi n+ih_n)$ is a vertical slit of height $2
h_n\ge 0, h_0=0$. The height $h_n$ is determined by the equation
$\cosh h_n=(-1)^n\D(e_n)\ge 1$, where $e_n\in [e_n^-,e_n^+]$ is such
that $\D'(e_n)=0$. The function $k$ maps the slit $g_n^c$ onto the
slit $\G_n$, and $k(-z)=-k(z)$ for all $z\in \cZ$.

In order to describe the spectral properties of the operator $H_0$, we start
from the properties of the canonical fundamental system $\vt, \vp$ of the
equation $-y''+py=z^2y$.
They satisfy the integral equations
\begin{multline}
\lb{fs}
\qqq\qqq\qqq\qqq  \vt(x,z)=\cos zx+\int_0^x{\sin z(x-s)\/z}p(s)\vt(s,z)ds,\\
\vp(x,z)={\sin zx\/z}+\int_0^x{\sin z(x-s)\/z}p(s)\vp(s,z)ds.\qqq\qqq\qqq\qqq
\end{multline}
For each fixed $x\in \R$, the functions $\vt(x,\cdot)$, $\vp(x,\cdot)$ are
entire in $\C$ and satisfy the estimates
\begin{multline}
\lb{efs1}
\max \rt\{|z|_1|\vp(x,z)|, \ |\vp'(x,z)| , |\vt(x,z)|,
{1\/|z|_1}|\vt'(x,z)|    \rt\} \le X=e^{|\Im z|x+\|p\|_x},\\
\rt|\vp(x,z)-{\sin zx\/z}\rt|\le {X\/|z|^2}\|p\|_x,
\qq |\vt(x,z)-{\cos zx}|\le {X\/|z|}\|p\|_x,
\end{multline}
where $|z|_1=\max\{1, |z|\}, \|p\|_x=\int_0^x|p(s)|ds$ and $(x,z)\in \R\ts \C$,
see page 13 in \cite{PT}. Substituting these estimates into \er{fs}  we obtain the standard
asymptotics
\[
\lb{asb} \b(z)={\vp'(1,z)-\vt(1,z)\/2}=\int_0^1{\sin z(2x-1)\/z}p(x)dx+{O(e^{|\Im z|})\/z^2}\qqq \as \qqq |z|\to \iy.
\]
Moreover, if $z=\pi n+O(1/n)$, then we obtain
\[
\lb{abn}
\b(z)=(-1)^n{p_{sn}\/2\pi n}+O({1/n^2}),\qq
p_{sn}=\int_0^1p(x)\sin 2\pi nx dx.
\]

The Floquet solutions $\p^{\pm}(x,z), z \in \cZ$, of the equation
$-y''+py=z^2y$ are given by
\[
\lb{3}
\p^\pm(x,z)=\vt(x,z)+m_\pm(z)\vp(x,z),\qqq \ m_\pm={\b\pm i\sin k\/ \vp(1,\cdot)},
\]
where $\vp(1,z)\p^+(\cdot,z)\in L^2(\R_+)$ for all
$z\in\C_+\cup\bigcup_{n\in\Z} g_n$.
In the trivial case $p=0$, we have $k=z$ and $\p^\pm(x,z)=e^{\pm izx}$.

Denote by $\cD_r(z_0)=\{|z-z_0|<r\}$ the disc of radius $r>0$ centred at
$z_0\in \L$. It is well known that if $g_n=\es$ for some $n\in \Z$, then the
functions $\sin k$ and $ m_\pm$ {\it are analytic in the disc}
$\cD_\ve(\m_n)\ss\cZ$ for some $\ve>0$, and the functions $\sin k$ and
$\vp(1,\cdot)$
have a simple zero at $\m_n$ \cite{F1}. Moreover, $m_\pm$ satisfies
\[
\lb{Tm-2}
m_\pm (\m_n)={\b'(\m_n)\pm i(-1)^nk'(\m_n)\/\pa_z\vp(1,\m_n)},
\qq \Im m_\pm (\m_n)\ne 0.
\]
Furthermore, $\Im m_+ (z)>0$ for all $z\in (z_{n-1}^+,z_{n}^-), n\in\N$,
and the following asymptotics hold,
\begin{equation}
\lb{Tm-1}
m_\pm (z)=\pm iz+O(1) \qq as \qq |z|\to \iy, \qq z\in \cZ_\ve, \ \ve > 0,
\end{equation}
where $\cZ_\ve =\{z\in \cZ\mid \dist \{z,g_n\}>\ve \ (n\in \Z, g_n\ne \es)\}$.
The function $\sin k$ and each function $\vp(1,\cdot)\p^{\pm}(x,\cdot)$,
$x\in \R$, are analytic on the Riemann surface $\cM$.
From \er{3}, the Floquet solutions
$\p^\pm(x,z)$ $((x,z)\in \R\ts \cM)$ satisfy \cite{T}
\[
\lb{f1}
\p^\pm(0,z)=1, \qq \p^\pm(0,z)'=m_\pm(z),
\qq \p^\pm(1,z)=e^{\pm ik(z)}, \qq  \p^\pm(1,z)'=e^{\pm ik(z)}m_\pm(z),
\]
\[
\lb{f2}
\p^\pm(x,z)=e^{\pm ik(z)x}(1+O(1/z))
\]
as $|z|\to \iy, z\in \cZ_\ve$, uniformly in $x\in\R$.
Below we shall require the simple identities
\[
\lb{LD0}
\b^2+1-\D^2=1-\vp'(1,\cdot)\vt(1,\cdot)= -\vp(1,\cdot)\vt'(1,\cdot).
\]

Now consider the modified Jost solutions $\P^\pm (x,z)$ of the equation
\[
\lb{bcf}
-{\P^\pm}''+(p+q)\P^\pm=z^2 \P^\pm \ {\rm on} \ [0,\infty)\qq {\rm with}\
\P^\pm(x,z)=\p^\pm(x,z), \ (x\ge t, z\in\cZ\sm \{0\}).
\]
Each function $\vp(1,\cdot)\P^\pm(x,\cdot)$, $x\ge 0$, is analytic in $\cM$,
since each $\vp(1,\cdot)\p^{\pm}(x,\cdot)$, $x\ge 0$, is  analytic in $\cM$.
We define the {\it modified Jost function} $\P_0^\pm=\P^\pm(0,z)$,
which is meromorphic in $\cM$ and has branch points $e_n^\pm$
$(n\in\Z, g_n\ne \es)$.

Let $\vp(x,z,\t)$ $((z,\t)\in \C\ts \R)$ be the solutions of the initial-value
problem
\[
\lb{x+t}
-\vp''+p(x+\t)\vp=z^2 \vp, \qq \ \vp(0,z,\t)=0,\qq \vp'(0,z,\t)=1.
\]

The solutions $\wt\vt, \wt\vp$ of the perturbed equation
$-y''+(p+q)y=z^2y, z\in \C$, determined by
$$
\wt\vt(x,z)=\vt(x,z),\ \wt\vp(x,z)=\vp(x,z)\ (\qq x\ge t),
$$
satisfy the integral equations
\begin{multline}
\lb{eqfwt}
\qqq\qqq\qqq\qq\wt\vt(x,z)=\vt(x,z)-\int_x^t\vp(x-s,z,s)q(s)\wt\vt(s,z)\,ds,\\
\wt\vp(x,z)=\vp(x,z)-\int_x^t\vp(x-s,z,s)q(s)\wt\vp(s,z)\,ds.\qqq\qqq\qqq\qqq
\end{multline}
For each $x\ge 0$, $\wt\vt(x,z)$, $\wt\vp(x,z)$ are entire functions of $z$ and
take real values if $z^2\in \R$.
The identities \er{eqfwt} and \er{bcf} give
\[
\lb{Ptp}
\P^\pm(x,z)=\wt\vt(x,z)+m_\pm(z)\wt\vp(x,z),\ ((x,z)\in \R_+\ts \cZ).
\]
We see that all singularities of $\P^\pm(x,z)$ coincide with some singularity
of $m_\pm(z)$ and hence do not depend on $x>0$.
The kernel of the resolvent $R=(H-z^2)^{-1},  z\in \C_+,$ then has the form
\[
\lb{R}
R(x,x',z)={\F (x,z )\P^+(x',z)\/\P_0^+(z)}\ ( x<x'),\ {\rm and}\
R(x,x',z )=R(x',x,z )\ (x>x').
\]
Here $\F(x,z)$ is the solution of the perturbed equation defined in \er{bcF}.
Each function $R(x,x',z),$ $ x,x'\in \R$, is meromorphic in $\cM$.
Note that if $z_0\in g_n^\pm\sm \{\m_n\pm i0\}$ for some $n$ and
$\P_0^+(z_0)\ne 0$, then the resolvent of $H$ is analytic at $z_0$.
The function $\P_0^+$ has a finite number of simple zeros
on each  $g_n^+, n\ne 0$, and on $i\R_+$, but no zeros on $\C_+\sm i\R_+$,
and the square of each zero is an eigenvalue. Note that the poles of
$$
\mR(x,z)={\P^+(x,z)\/\P^+(0,z)}={\wt \vt(x,z)+m_+(z)\wt \vt(x,z)\/\wt \vt(0,z)+m_+(z)\wt \vt(0,z)}
$$
are either poles of $m_+$ or zeros of $\P_0^+$ and therefore are independent of
$x$.
A pole of $\mR(x,\cdot)=\P^+(x,\cdot)/\P_0^+$ on $g_n^+$ is called a {\it bound
state} by slight abuse of terminology (since it is a momentum rather than an
energy value). Moreover, $\P_0^+(z)$ has an infinite number of zeros in  $\C_-$, see \er{T2-3}.
In terms of the Riemann surface $\cM$, the previous definition of resonances
and states, which referred to the Riemann surface $\L$, takes the following
form.

\bigskip
\no {\bf Definition. }
1) A point $\z\in \ol\C_-\cap \cM, \z\ne 0$, is a {\it resonance} if the
function $\mR(x,z)$ has a pole at $\z$  for almost all $x>0$.

\no 2) A point  $\z=e_n^\pm , n\ne 0$, (or $\z=0$) is a {\it virtual state}
if the function $z\mapsto\mR(x,\z+z^2)$ (or $\mR(x,z)$) has a pole at $0$ for
almost all $x>0$.

\no 3) A point $\z\in\cM$ is a {\it state} if it is either a bound state or
a resonance or a virtual state.
If $\z\in g_n^-, n\ne 0$, or $\z\in g_0^-=i\R_-$, then we call $\z$ an
{\it antibound state}.

\bigskip
Let  $\s_{bs}(H)$ (or $\s_{rs}(H)$ or $\s_{vs}(H)$ ) be the set of all bound
states (or resonances or virtual states) of $H$ in $\cM$ (i.e., in terms of the
momentum variable), and set $\s_{st}(H)=\s_{vs}(H)\cup\s_{rs}(H)\cup\s_{bs}(H)$.

The kernel of the resolvent $R_0(z)=(H_0-z^2)^{-1}, z\in \C_+$, has the form
\[
\lb{R0}
R_0(x,x',z)=\vp(x,z)\p^+(x',z)\ (x<x'),\ \
{\rm and} \qq R_0(x,x',z)=R_0(x',x,z)\ (x>x').
\]
Consider the states $z_n^0=\sqrt{\l_n^0}\in g_n^c$ of $H_0$.
Due to \er{3}, the function $\P_0^+=\p^+(0,\cdot)=1$ and
$\mR(x,\cdot)=\vt(x,\cdot)+m_+\vp(x,\cdot)$.
Since $\vp(1,\cdot)m_+$ and $\sin k(z)$ are analytic in $\cM$,
the resolvent $R_0(z)$ has singularities only at $\m_n\pm i0$, where $g_n\ne \es$,
and in order to describe the states of $H_0$ we need to study $m_+$ on $g_n^c$ only.
We need the following result (see \cite{Zh3}),
where $\mA(z_0), z_0\in \cM$, denotes the set of functions analytic in some
disc $\cD_r(z_0), r>0$.

\begin{lemma}
\lb{Tm}
All states of $H_0$ are of the form $\m_n\pm i0\in g_n^c, n\ne 0$, where
$g_n\ne \es$.
Consider a non-empty momentum gap $g_n=(e_n^-,e_n^+)$ for some $n\ge 1$.
Then

\no i)  $z_n^0=\m_n+i0\in g_n^+$ is a bound state of $H_0$ if and only if one
of the following conditions is true.
\begin{multline}
\lb{Tm-31}
(1)\qq  m_-\in \mA(\m_n+i0),\\
\ \ (2)\qq \b(\m_n)=i\sin k(\m_n+i0)= -(-1)^n\sinh h_{sn},  \qq  k(\m_n+i0)=\pi n+ih_{sn}\qq h_{sn}>0,  \hspace{1.2cm} \\
(3)\qq  m_+(z_n^0+z)={c_n\/z}+O(1) \qq as \ z\to 0, \ z\in \C_+,\
c_n={-2\sinh |h_{sn}|\/(-1)^n\pa_z\vp(1,\m_n)}<0.\qqq
\end{multline}
\no ii)  $z_n^0=\m_n-i0\in g_n^-$ is an antibound state of $H_0$ if and only if
one of the following conditions is true.
\begin{multline}
\lb{Tm-32}
(1) \qqq   m_-\in \mA(\m_n-i0),\\
(2)\qqq
\b(\m_n)=i\sin k(\m_n-i0)= -(-1)^n\sinh h_{sn},  \qq k(\m_n-i0)=\pi n+ih_{sn},\qq h_{sn}<0,    \hspace{0.3cm} \\
(3)\qqq  \qq m_+(z_n^0+z)={-c_n\/z}+O(1) \qq as \ z\to 0, \ \ z\in \C_-.  \hspace{4.5cm}
\end{multline}

\no iii) $z_n^0=\m_n$ is a virtual state of $H_0$ if and only if one of the
following conditions is true.
\begin{multline}
\lb{Tm-33}
(1) \qqq z_n^0=\m_n=e_n^-\qqq  or \qqq z_n^0=\m_n=e_n^+,\\
(2) \qqq m_+(z_n^0+z)={c_n^0\/\sqrt z}+O(\sqrt z)\qq as \ z\to 0, \ z\in \C_+,\qq
c_n^0\ne 0. \qqq\qqq\qqq\qq\
\end{multline}
\end{lemma}

\noindent
These facts are well known in inverse spectral theory (\cite{N-Z},
\cite{MO}, \cite{K5}).
A detailed analysis of $H_0$ was done in \cite{Zh3}, \cite{K4}.

If $\m_n\in g_n\ne \es$, then the function $m_+$ has a pole at
$z_n^0=\m_n+i0\in g_n^+$ (corresponding to a bound state) or at
$z_n^0=\m_n-i0\in g_n^-$ (corresponding to an antibound state).
If $\m_n=e_n^+$ (or $\m_n=e_n^-$), then $z_n^0=\m_n$ is a virtual state.
Note that for closed gaps $g_n=\es, n\ne 0$, each $\p^\pm(x,\cdot)\in
\mA(\m_n), x\ge 0$.
Moreover, the  resolvent $R_0(z)$ has a pole at $z_0$ if and only if the
function $m_+(\cdot)$ has a pole at $z_0$.

 The following asymptotics from \cite{MO} hold true as $n\to \iy$:
\[
\lb{sde}
\m_n=\pi n+\ve_n(p_{0}-p_{cn}+O(\ve_n)),\qqq p_{cn}=\int_0^1p(x)\cos 2\pi nx\,dx,
\qqq
\qq \ve_n={1\/2\pi n},
\]
\[
\lb{anc}
h_{sn}=-\ve_n(p_{sn}+O(\ve_n)),
\]
\[
\lb{ape}
e_n^\pm=\pi n+\ve_n(p_0\pm |p_n|+O(\ve_n)), \qq \qq
p_n=\int_0^1p(x)e^{-i2\pi nx}dx=p_{cn}-ip_{sn}.
\]
The fundamental system $y_1$, $y_2$ of the equation $-y''+(p+q)y=z^2y, z\in \C$,
with initial data
\[
\lb{wtc}
y_2'(t,z)=y_1(t,z)=1,\qqq y_2(t,z)=y_1'(t,z)=0,
\]
satisfies the integral equations
\begin{multline}
\lb{efy}
\qqq\qqq
y_1(x,z)=\cos z(x-t)-\int_x^t{\sin z(x-\t)\/z}(p(\t)+q(\t))y_1(\t,z)d\t,\\
y_2(x,z)={\sin z(x-t)\/z}-\int_x^t{\sin z(x-\t)\/z}(p(\t)+q(\t))y_2(\t,z)d\t.
\qqq\qqq\
\end{multline}
For each $x\in \R$ the functions  $ y_1(x,z),  y_2(x,z)$  are entire in $z\in\C$ and satisfy
\begin{multline}
\lb{efs}
\max \{||z|_1 y_2(x,z)|, \ | y_2'(x,z)| , | y_1(x,z)|,
{1\/|z|_1}| y_1'(x,z)|    \} \le X_1=e^{|\Im z||t-x|+\int_x^t|p+q|ds},\\
| y_1(x,z)-\cos z(x-t)|\le {X_1\/|z|}\|q\|_t,\qq
\rt| y_2(x,z)-{\sin z(x-t)\/z}\rt|\le {X_1\/|z|^2}\|q\|_t,
\end{multline}
in analogy to \er{efs1}.

The initial-value problem for the inhomogeneous equation $-y''+(p-z^2)y=f,\
 y(0)=y'(0)=0$ has the unique solution
$
y(x)=-\int_0^x\vp(x-\t,z,\t)f(\t)d\t.
$
Hence the solutions $\F$ and $\P^\pm$ of $-y''+(p+q)y=z^2y$ satisfy
\[
\lb{eF}
\F(x,z)=\vp(x,z)+\int_0^x\vp(x-s,z,s)q(s)\F (s,z)\,ds,
\]
\[
\lb{ep}
\P^\pm (x,z)=\p^\pm(x,z)-\int_x^t\vp(x-s,z,s)q(s)\P^\pm(s,z)\,ds.
\]
Below we need the well known fact from scattering theory that
\[
\lb{DP}
\P^+(0,z)=D(z^2)=\det (I+q(H_0-z^2)^{-1})\ (z\in \cM).
\]
This is similar to the case $p=0$, see \cite{J}. The case on the
real line with $p\ne \const$ was considered in \cite{F4}. The
functions $k, \P^\pm, m_\pm, \p^\pm$ are meromorphic in $\cZ$ and
real on $i\R$ and then they satisfy
\[
k(-z)=\ol k(\ol z), \qq \P^\pm(\cdot,-z)=\ol \P^\pm(\cdot,\ol z),
\qq m_\pm(-z)=\ol m_\pm(\ol z), \qq \p^\pm(\cdot,-z)=\ol
\p^\pm(\cdot,\ol z)
\]
for all $z\in \cZ$.

\begin{lemma}
\lb{T21}
i) The following identities, estimate and asymptotics  hold true.
\[
\lb{T21-1}
\P^\pm =\p^{\pm}(t,\cdot)y_1+\pa_t\p^{\pm}(t,\cdot)y_2,
\]
\[
\lb{T21-2}
\P^\pm(0,z)=1+\int_0^t\vp(x,z)q(x)\P^\pm(x,z)\,dx,
\]
\[
\lb{T21-200}
|\P^\pm(x,z)-\p^\pm(x,z)|\le e^{\mp vx+b(t,x)}{w^\pm(z)\/|z|_1}\int_x^t|q(r)|
\,dr,\qq \ (x\in [0,t]),
\]
where
$$
v=\Im k(z), \qq w^\pm=\sup_{x\in [0,t]} |e^{\pm vx}\p^\pm(x,z)|, \qq b(t,x)=(|v|-v)(t-x)+\int_x^t (|p(r)|+|q(r)|)\,dr
$$
and
\[
\lb{T21-3}
\P^\pm(x,z)=e^{\pm ik(z)x}(1+e^{\pm (t-x)(|v|-v)}O(1/z))
\]
as $|z|\to \iy, z\in \cZ_\ve$, $\ve>0$, uniformly in $x\in [0,t]$.
Moreover, {\rm \er{T1-1}} holds true.

\no ii) The function $\P^\pm(0,\cdot)$ has
exponential type $2t$ in the half plane $\C_\mp$.

\end{lemma}

\noindent
An entire function $f(z)$ is said to be {\it of exponential type} if there is
a constant $A$ such that  $|f(z)|\le $ const $e^{A|z|}$ everywhere [Koo].
The infimum over the set of $A$ for which such an inequality holds is called
the {\it type} of $f$.

\medskip
\no{\bf Proof.} i)  Using \er{wtc}, \er{bcf} we obtain \er{T21-1}.

The identity $\vp(x,z,\t)=\vt(\t,z)\vp(x+\t,z)-\vp(\t,z)\vt(x+\t,z)$ gives
$\vp(-x,z,x)=-\vp(x,z)$ and \er{ep} yield \er{T21-2}.

The estimate \er{T21-200} was proved in \cite{K4}.
Substituting \er{f2} into \er{T21-200}  we obtain \er{T21-3}.

In particular, substitution of \er{T21-3}, \er{efs} into \er{T21-2} yields
\er{T1-1}, since $D(z^2)=\P^+(0,z)$.

ii) We give the proof for the case $\P^+(0,z)$, the proof for $\P^-(0,z)$ being
similar. Due to \er{T21-3}, $\P^+(0,z)$ has exponential type $\le 2t$ in the half plane $\C_-$.
The decompositions
$$
f(x,z)\ev e^{-ixk(z)}\P^+(x,z)=1+\r f_1(x,z), \qqq
\vp(x,z)e^{ixk(z)}\ev \r (e^{i2xz}-1 +\r\e(x,z)), \
$$
where $\r={1\/2iz}$,  give
\begin{multline}
\lb{esff}
\P^+(0,z)-1=\int_0^t\vp(x,z)e^{ixk}q(x)f(x,z)\,dx\\
\qq =
\r\int_0^te^{i2xz}q(x)f(x,z)(1+\r e^{-i2xz}\e(x,z))\,dx-\r\int_0^tq(x)f(x,z)\,dx\\
=\r K(z)-\r\int_0^tq(x)\,dx,\qqq\qqq\qqq\qqq\qqq\qqq\qqq\qqq\qq\\
K=\int_0^te^{i2xz}q(x)(1+G(x,z))\,dx, \qqq G=\r f_1+\r e^{-i2xz}(\e
f-f_1), \qq  \  z\in \cZ_\ve.
\end{multline}

Asymptotics \er{efs1}, \er{T21-3} and $k(z)=z+O(1/z)$ as $z\to \iy$ (see \cite{KK}) yield
\[
\lb{esff1}
\e(x,z)= e^{2x|\Im z|}O(1),\qq \qqq f_1(x,z)=e^{2(t-x)|\Im
z|}O(1)\qq as \ |z|\to \iy, \ z\in \cZ_\ve.
\]

Now we use the following variant of the Paley-Wiener Theorem from
\cite{Fr}.

\no {\it Let $h\in \cQ_t$ and let each $F(x,z), x\in [0,t]$, be analytic for
$z\in\C_-$ and $F\in L^2((0,t)\times\R)$.
Then $\int_0^te^{2izx}h(x)(1+F(x,z))\,dx$ has exponential type at least $2t$ in
$\C_-$.}

We cannot apply this result to the function  $K(z), z\in\C_-$, since
$\p^+(x,z)$ has a singularity at $z_n^0$ if $g_n\ne \es$.
However, we can use it for the function  $K(z-i), z\in\C_-$,
since \er{esff}, \er{esff1} imply that $\sup _{x\in [0,1]}|G(x,-i+\t)|=O(1/\t)$
as $\t\to \pm\iy$.
Thus the function $\P^+(0,z)$ has exponential type $2t$ in the half plane
$\C_-$.
 \BBox

\section {Spectral properties of $H$}
\setcounter{equation}{0}

For the sake of brevity, we write
$$ \vp_1=\vp(1,\cdot ),\ \vp_1'=\vp'(1,\cdot ), \ \vt_1=\vt(1,\cdot ),
\ \vt_1'=\vt'(1,\cdot ), \ \F_1=\F(1,\cdot), \ \F_1'=\F'(1,\cdot),
\ \psi_t^\pm = \psi^\pm(t,\cdot).
$$
We define the function
\[
F(z)=\vp(1,z)\P^-(0,z)\P^+(0,z) \ (z\in \cZ),
\]
which is real on $\R$, since $\P^-(0,z)=\ol\P^+(0,\ol z)$ for all $z\in \cZ$,
see also \er{T22-1}.

A function $f$ is said to belong to the {\it Cartwright class} $Cart_\o$
if $f$ is entire, of exponential type and satisfies:
$$\o_\pm(f):=\limsup_{y\to \iy}
{\log |f(\pm iy)|\/y}=\o>0,\qqq \int_{\R}{\log (1+|f(x)|)dx\/ 1+x^2}<\iy.
$$

In the following, the partial derivative w.r.t. $t$ is denoted by a dot.

\begin{lemma} \lb{T22}
i) The following identities and estimates hold true.
\[
\lb{T22-1}
F=\vp(1,\cdot,t)y_1^2(0,\cdot)+\dot\vp(1,\cdot,t)y_1(0,\cdot)y_2(0,\cdot)-\vt'(1,\cdot,t)y_2^2(0,\cdot)\in Cart_{1+2t},
\]
\[
\lb{T22-2}
\vp_1\dot\p_t^+\dot\p_t^-=-\vt'(1,\cdot,t),\qqq
\vp_1(\dot\p_t^+\p_t^-+\p_t^+\dot\p_t^-)=\dot\vp(1,\cdot,t)=
\vp'(1,\cdot,t)-\vt(1,\cdot,t),
\]
\[
\lb{T22-3}
|F(z)-{\sin z\/z}|\le {C_F e^{(1+2t)|\Im z|}\/|z|^2},
\qq C_F=3(\|p\|_1+\|p+q\|_t)e^{2\|p+q\|_t+\|p\|_1}.
\]

\no ii) The set of zeros of $F$ is
symmetric w.\ r.\ t.\ both the real line and the imaginary line.
In each disc $\{z\mid |z-\pi n|<{\pi\/ 4}\}, |n|\ge 1+4C_Fe^{t{\pi\/2}}$, there
exists exactly one simple real zero $z_n$ of $F$, and $F$ has no zeros in the
domain $\mD_F\cap\C_-$, where $\mD_F=\{z\in \C\mid 4C_Fe^{2|\Im z|}<|z|\}$.

\no iii) For all $z \in \cZ$,
 \[
\lb{T31-1}
\P_0^\pm(z)=e^{\pm ik(z)n_t}w_\pm(z),\ \ \  w_\pm(z)=\F'(n_t,z)-m_\pm(z)\F(n_t,z).
\]

\end{lemma}
\no{\bf Proof.} i)
The function $\vp(1,z,t)$, $(t,z)\in \R\ts \C$, (see \er{x+t}) satisfies
\cite{Tr}
\[
\lb{ipv}
\vp(1,\cdot,t)=-\vt_1'\vp_t^2+\vp_1\vt_t^2+
2\b\vp_t\vt_t=\vp_1\p_t^+\p_t^-.
\]
Hence  we obtain
\begin{multline}
\ \qqq
\lb{2is}
\dot \vp(1,\cdot,t)=\vp_1(\dot \p_t^+\p_t^-+\p_t^+\dot\p_t^-),\\
\ddot \vp(1,\cdot,t)=\vp_1(\ddot \p_t^+\p_t^-+\p_t^+\ddot\p_t^-
+2\dot \p_t^+\dot \p_t^-)=2\vp_1(p(t)-z^2)\p_t^+\p_t^-+
2\vp_1\dot \p_t^+\dot \p_t^-.
\ \qq
\end{multline}
Identity \er{T21-1} gives
\[
\lb{e12}
F=\vp_1(\p_t^+\p_t^-y_1^2(0,\cdot)+\dot\p_t^+\dot\p_t^-y_2^2(0,\cdot)+
(\p_t^+\dot\p_t^-+\dot\p_t^+\p_t^-)y_1(0,\cdot)y_2(0,\cdot).
\]
Then using the following identities from \cite{IM},
\[
\lb{im1}
\dot \vp(1,z,t)=\vp'(1,z,t)-\vt(1,z,t),\qq
\ddot \vp(1,z,t)=2(p(t)-z^2)\vp(1,z,t)-2\vt'(1,z,t),
\]
and  \er{ipv}, \er{e12}, we obtain \er{T22-1}, \er{T22-2}, since Lemma
\ref{T21}, ii) and \er{efs1} yield $F\in Cart_{1+2t}$.

Next we show \er{T22-3}. We have
\[
y_1(0,\cdot)=\cos t z+\wt y_1,\qqq \wt y_1=\int_0^t
{\sin z s\/z}(p(s)+q(s))y_1(s,\cdot)\,ds,
\]
\[
y_2(0,\cdot)=-{\sin t z\/z}+\wt y_2,\qqq \wt y_2=\int_0^t
{\sin z s\/z}(p(s)+q(s))y_2(s,\cdot)\,ds,
\]
\[
\vt(1,t)=\cos z+\vt_{1t},\qqq \vt_{1t}=\int_0^1
{\sin z(1-s)\/z}p(s+t)\vt(s,t)\,ds,
\]
\[
\vt'(1,t)=-z\sin z+\vt_{1t}',\qqq \vt_{1t}'=\int_0^1
\cos z(1-s)p(s+t)\vt(s,t)\,ds,
\]
\[
\vp(1,t)={\sin z\/z}+\vp_{1t},\qqq \vp_{1t}=\int_0^1
{\sin z(1-s)\/z}p(s+t)\vp(s,t)\,ds.
\]
Then \er{T22-1} implies
\begin{multline}
F=(\cos tz+\wt y_1)^2({\sin z\/z}+\vp_{1t})+({\sin tz\/z}-\wt y_2)^2(z\sin z-\vt_{1t}')
+(\cos tz+\wt y_1)(-{\sin tz\/z}+\wt y_2)\dot \vp(1,z,t)\\
={\sin t z\/z}+f_1+f_2+f_3,
\qqq\qqq\qqq\qqq\qqq\qqq\qqq\qqq\qqq\qqq\qqq\qqq\qq
\end{multline}
where
$$
 f_1=y_1(0,\cdot)^2\vp_{1t}+{\sin tz\/z}(\cos tz+y_1(0,\cdot) )\wt y_1,
$$$$
f_2=-y_2(0,\cdot)^2\vt_{1t}'+z\sin z(y_{2}(0,\cdot)-{\sin tz\/z})\wt y_2,
\qqq  f_3=y_1(0,\cdot)y_2(0,\cdot)\dot \vp(1,z,t).
$$
The estimates
$$
\qq
 |f_3|\le {C_t\/|z|^2}\|p\|_1,\qq
|f_j|\le {C_t\/|z|^2}(\|p\|_1+2\|p+q\|_t)\ (j\in\{1,2\})
$$
yield \er{T22-3}, where
$C_t=e^{(2t+1)|\Im z|+2\|p+q\|_t+\|p\|_1}$.

ii) From \er{T22-3} we obtain for $|n|\ge 1+4C_Fe^{t{\pi\/2}}$
$$
|F(z)-{\sin z\/z}|\le {C_F\/|z|^2}e^{|\Im z|+t{\pi\/2}|}
\le {4C_F\/|z|}e^{t{\pi\/2}|}{|\sin z|\/|z|}<{|\sin z|\/|z|}
$$
for all $|z-\pi n|={\pi\/4}$,
since $e^{|{\Im}z|}\le 4|\sin z|$ for all
$|z-\pi n|\ge \pi /4, n\in \Z$ (see [PT]).
Hence, by Rouch\'e's theorem, $F$ has the same number of roots, counting
multiplicities, as $\sin z$ in each disc $\cD_{\pi\/4}(\pi n)$.
Since $\sin z$ only has the roots $\pi n$ $(n\in\Z)$, the statement about the
zeros of $F$ in each disc follows. The zero in $\cD_{\pi\/4}(\pi n)$ is real,
since $F$ is real-valued on the real line.

Using \er{T22-3} and $e^{|{\Im}z|}\le 4|\sin z|$ for all
$|z-\pi n|\ge \pi /4, n\in \Z$, we obtain
$$
|F(z)|\ge |{\sin z\/z}|-\rt|F(z)-{\sin z\/z}\rt|\ge {e^{|\Im z|}\/4|z|^2}
\rt(|z|-4C_F e^{2t|\Im z|}   \rt)>0
$$
for all $z\in \mD_1=\{z\in \mD_Fz\mid |z-\pi n|\ge \pi /4 (n\in \Z) \}$.
As the function $F$ has exactly one real zero $z_n$ in $\cD_{\pi\/4}(\pi n),
n\ge n_0$, it therefore has no zeros in the domain $\mD_F$.
That  the set of zeros of $F$ is symmetric with respect to both the real line
and the imaginary line follows from the fact that $F$ is real-valued on both
lines.

iii) Using $\P^\pm(0,z)=\{\P^\pm(x,z ),\F (x,z)\}$ $(z \in \cZ)$ at $x=n_t$,
and  \er{f1} we obtain  \er{T31-1}, where   $\{y,u\}=yu'-y'u$ is the Wronskian.
\BBox

\medskip

With $g_n^c={\ol g_n^-}\cup {\ol g_n^+}$ as before, we define the sets
\[
\mZ=i\R\cup \C_-\cup \bigcup_{n\in \Z}g_n^c \subset\cZ,\qqq
\mZ_0=\mZ\sm \{0, e_n^\pm,\m_n\pm i0 (n\in \Z, g_n\ne \es)\}.
\]
We describe the zeros of $\P^\pm (0,z)$.

\begin{lemma} \lb{T311}
i) If $g_n=(e_n^-,e_n^+)= \es$ for some $n\ne 0$, then $\P^\pm(x,\cdot )\in \mA(\m_n)$ $(x\ge 0)$, and
$\Im \P^\pm(0,\m_n)\ne 0$. Moreover, $\m_n=e_n^\pm$ is a simple zero of $F$
and $\m_n\notin\s_{st}(H)$.

\no ii)  $\P^\pm(0,z)\ne 0$ for all $z\in (e_{n-1}^+,e_{n}^-), n\in\Z$.
Moreover,  states of $H$ and zeros of $\P^+(0,z)$ lie in $\mZ$.

\no iii) A point $z\in\mZ_0$
  is a zero of $\P_0^+$ if and only if $z\in\mZ_0\cap \gS_{st}(H)$.
In particular, if $z\in\mZ_0\cap \gS_{st}(H)$, then

  1) $z\in \ol\C_+\cap \cZ$ is a bound state of $H$,

   2) $z\in \ol\C_-\cap \cZ$  is a resonance of $H$.

\end{lemma}
\no {\bf Proof.} i)
Lemma \ref{Tm} and identity \er{Ptp} yield that $\P^\pm(x,\cdot )\in \mA(\m_n)$
$(x\ge 0)$.

Using \er{Tm-2} we deduce that $\P_0^\pm(\m_n)\ne 0$,
hence $\P^\pm(x,\cdot), 1/\P_0^\pm\in \mA(\m_n)$. Thus $\m_n$ is not a  state of $H$ and $\m_n$ is a simple zero of $F$.

 ii) The conformal mapping $k$ maps each interval $(e_{n-1}^+,e_{n}^-),
n\ge 1$, onto $(\pi (n-1), \pi n)$.
As $\m_n\in [e_{n}^-,e_{n}^+]$, the function $m_\pm$ is analytic on
$(e_{n-1}^+,e_{n}^-)$ and $\Im m_\pm(z)\ne 0$ for all
$z\in (e_{n-1}^+,e_{n}^-)$. Then the identity \er{Ptp}
gives $\P^\pm(0,z)\ne 0$ for all $z\in (e_{n-1}^+,e_{n}^-)$.

The identity $\P^+(0,z)=D(z^2)$ and standard arguments (similar to the case
$p=0$, see \cite{K1})
imply that states of $H$ and zeros of $\P^+(0,z)$  belong to the set
$\mZ$.

 iii) follows from the identities  \er{Ptp},  \er{R}.
$\BBox$

\medskip
\noindent
States of $H$ which coincide with unperturbed states $z_n^0$ have the following
properties.

\begin{lemma}
\lb{T32}
Let $\z=\m_n+i0\in g_n^+$ or $\z=\m_n-i0\in g_n^-$  for some $n\ge 1$,
where $g_n\ne \es$. Then

\no i) $\wt \vp(0,\m_n)=0$ if and only if $\F(n_t,\m_n)=0$.

\no ii) Assume in addition $\z=z_n^0\in \s_{st}(H_0)$.
Then $\P_0^-\in \mA(\z)$, each $\P^+(x,\cdot)$, $x>0$, has a simple pole at
$\z$, and there are two cases:

\no 1) if $\wt\vp(0,\m_n)=0$, then $\P_0^+\in\mA(\z), \ \P_0^-(\z)\ne 0$
and $\z\in \s_{st}(H)$. In particular,
\[
\lb{T32-0}
\qq \z=\m_n+i0\in g_n^+ \qq \Rightarrow \qq \P_0^+(\z)\ne 0,\qq
F(\m_n)=0,\qqq (-1)^nF'(\m_n)>0;
\]
2) if $\wt\vp(0,\m_n)\ne 0$,  then
\begin{multline}
\lb{T32-1}
\qqq\qqq\qqq
{\P^+(x,\cdot)\/\P_0^+}\in \mA(\z) \ (x\ge 0), \ \  \qq
\z\notin \gS_{st}(H),\qq \P_0^-(\z)\ne 0,\qq
F(\m_n)\ne 0,\\
\P_0^+(z)={c_n+O(\ve)\/\ve}\wt\vp(0)\qq as \ \ve=z-\z\to 0.
\qqq\qqq\qqq
\end{multline}

\no iii) $\z\in \gS_{bs}(H) $ (or $\z\in \s_{vs}(H) $ ) if and only if
$\z\in \gS_{bs}(H_0) $ (or $\z\in \s_{vs}(H_0)$) and $\F(n_t,\m_n)=0$.

\no iv) If $\z\in \gS_{st}(H_0)\cap \gS_{st}(H)$, then
its complex conjugate $\ol \z$ on the other sheet is not a state of $H$.

\no v) Let $\z\in \gS_{st}(H_0)$ and assume that its complex conjugate $\ol\z$
on the other sheet is a state of $H$. Then $\z\notin \gS_{st}(H)$.

\end{lemma}
\no{\bf Proof.} i) Comparing \er{T31-1} and \er{Ptp} we deduce that
 $\wt \vp(0,\m_n)=0$ if and only if $\F(0,\m_n)=0$.

ii) Lemma \ref{Tm} yields $m_-\in \mA(\z)$
and each $\P^+(x,\cdot)$, $x>0$, has a simple pole at $\z$
and $m_+(z)={\pm c_n\/\ve}+O(1)$ as $\ve\to 0,\ \ve \in \C_\pm$, with $c_n<0$.

1) If $\wt\vp(0,\m_n)=0$, then \er{Ptp} yields     $\P_0^+\in\mA(\z), \ \P_0^-(\z)\ne 0$.
Thus $\z\in \gS_{st}(H)$, since each $\P^+(x,\cdot)$, $x>0$< has a simple pole
at $\z$.

Consider the case $\z=\m_n+i0\in g_n^+$. As \er{T31-1} gives
$$
\P^\pm(0,z)=e^{\pm ik(z)n_t}w_\pm(z),\ \ \  w_\pm(z)=\F'(n_t,z)-m_\pm(z)\F(n_t,z),
\qqq (z\in \cZ),
$$
we find
$$
w_-(\z)=\F'(n_t,\m_n)\ne 0,\qq
w_+(z_1)=\F'(n_t,\m_n)(1-c_n{\pa_z\F(n_t,\m_n)\/\F'(n_t,\m_n)})\ne 0,
$$
since  $\F'(n_t,\m_n)\ne 0, \F(n_t,\m_n)=0$  and $\F'(n_t,\m_n)\pa_z\F(n_t,\m_n)>0$  \cite{PT}.
This yields \er{T32-0}, since $(-1)^n\pa_z \vp(1,\m_n)>0$
\cite{PT}.

2) Lemma \ref{Tm} gives $m_+(z)={c_n+O(\ve)\/\ve}$ as $\ve=z-\z\to 0$.
Using \er{Ptp}, we obtain
$$
\P_0^+(z)={c_n+O(\ve)\/\ve}\wt\vp(0),\qq
{\P^+(x,z)\/\P_0^+(z)}={\ve \wt\vt(x)-c_n\wt\vp(x)+O(\ve)\/\ve \wt\vt(0)-c_n\wt\vp(0)+O(\ve)}=
{\wt\vp(x)\/\wt\vp(0)}+O(\ve)\qq \as \ \ve\to 0,
$$
since $c_n\wt\vp(0,\m_n)\ne 0$, where we abbreviate $\wt\vp(x)=\wt\vp(x,z)$ etc.
This yields $\vp(1,z)\P_0^+(z)=\pa_z\vp(1,\m_n)c_n+o(1)$, and
$m_-\in \mA(\z)$ gives $\P_0^+(\z)=\wt\vt(0)\ne 0$, so $F(\m_n)\ne 0$
and \er{T32-1}.

Using   i) and  ii) we obtain iii).

\no iv) $\P_0^+\in \mA(\ol\z)$
and each $\P^+(x,\cdot)\in \mA(\ol\z), x>0$.
Due to ii) $\wt\vp(0,\m_m)=0$, then we obtain $\P_0^+(\ol\z)\ne 0$.
Thus $\ol\z\notin \gS_{st}(H)$.

\no v) follows directly from iv) by contraposition.
\BBox

\medskip
\noindent
For virtual states, which coincide with the points $e_n^\pm$, we have the
following.

\begin{lemma} \lb{T33}
Let $\z=e_n^-$ or $\z=e_n^+$ for some $n\ge 1$, where $e_n^-<e_n^+$ and set $\ve=z-\z$.

\no i) Assume that $\z\ne \m_n$ and $\P_0^+(\z)=0$. Then
 $\z$ is a simple zero of $F$, $\z\in \s_{vs}(H)$
 and
\[
\lb{T33-1}
\P_0^+(z)=\wt\vp(0,\z)c\sqrt \ve+O(\ve),\qq
\mR(x,z)={\P^+(x,z)+O(\ve)\/\wt\vp(0,\z)c\sqrt \ve},\qq c\wt\vp(0,\z)\ne 0.
\]
\no ii) Assume that $\z=\m_n$ and $\wt\vp(0,\z)\ne 0$. Then
$F(\z)\ne 0$ and each $\mR(x,\cdot), x>0$, does not have a singularity at $\z$,
and $\z\notin \s_{vs}(H)$.

\no iii) Assume that $\z=\m_n$ and $\wt\vp(0,\z)=0$. Then $\z\in \s_{vs}(H)$,
$\P_0^\pm(\z)\ne 0$ and $\z$ is a simple zero of $F$, and each
$\mR^2(x,\cdot), x>0$, has  a pole at $\z$.
\end{lemma}

\no {\bf Proof.}
i) Lemma \ref{Tm} gives $m_\pm(z)=m_\pm(\z)+c\sqrt \ve+O(\ve)$ as $\ve=z-\z\to 0, c\ne 0$.
We distinguish two cases. First assume $\wt\vp(0,\z)\ne 0$.
Then the identity \er{Ptp} implies \er{T33-1}.

Second, if $\wt\vp(0,\z)=0$, then \er{Ptp} implies
 $\P_0^+(\z)=\wt\vt(0,\z)\ne 0$, which gives a contradiction.

\no ii) If $\z=\m_n$, then Lemma \ref{Tm} gives
$m_\pm(z)=\pm {c\/\sqrt \ve}+O(1)$ $(\ve\to 0)$ with $c\ne 0$.
Then  \er{Ptp} implies
$$
\P_0^\pm(z)=\pm{\wt\vp(0,\z)c\/\sqrt \ve}+O(1),\qq
{\P^+(x,z)\/\P_0^+(z)}={\wt\vt(x,z)+({c\/\sqrt \ve}+O(1)) \wt\vp(x,z)\/{\wt\vp(0,\z)c\/\sqrt \ve}+O(1)}={1+O(\sqrt \ve)\/\wt\vp(0,\z)}.
$$
Thus the function $\mR(x,\cdot), x>0$, is not singular at $\z$, and
$\z\notin \s_{vs}(H)$, $F(\z)\ne 0$.

iii) If $\wt\vp(0,\z)=0$, then \er{Ptp} gives $\P_0^+(\z)=\wt\vt(0,\z)\ne 0$, since $\wt\vt(0,\z)\ne 0$ and $\b(\z)=0$. Moreover, we obtain
$
\P^+(x,z)=\wt\vt(x,z)+({c\/\sqrt \ve}+O(1)) \wt\vp(x,z),
$
and the function $\mR^2(x,\cdot), x>0$, has a pole at $\z$, $\z\in \s_{vs}(H)$
and $F(\z)=0$.
$\BBox$

\begin{lemma}
\lb{T34}
Let $\l\in \g_n$, $\l\ne \m_n^2$, be an eigenvalue of $H$ for some $n\ge 0$ and
let $z=\sqrt \l\in i\R_+\cup \bigcup_{n\ge 1} g_n^+$. Then
\[
\lb{T34-1}
C_\l=\int_0^\iy |\P^+(x,z)|^2\,dx=-{{\P^+}'(0,z)\/2z} \pa_z\P^+(0,z)>0,
\]
\[
\lb{T34-2}
{2i\sin k(z)\/\vp(1,z)}=\P^-(0,z_1){\P^+}'(0,z)\ne 0,\qqq
i\sin k(z)=-(-1)^n\sinh h,\ \ h=h(\l)>0,
\]
\[
\lb{T34-3}
C_\l={(-1)^nF'(z)\sinh h\/ z\vp^2(1,z)\P^-(0,z)^2}>0,
\qqq {(-1)^nF'(z)\/z}>0.
\]
\end{lemma}
\no {\bf Proof.} \er{T34-1} follows from the identity
$\{{\pa\/\pa z} \P^+,\P^+\}'=2z(\P^+)^2$.

Using the Wronskian of the functions $\P^+, \P^-$ and \er{f1}, we obtain
$\P^-(0,z){\P^+}'(0,z)=m_+(z)-m_-(z)$,
which yields \er{T34-2}, since $k(z)=\pi n+ih$ for some $h>0$ (see the
definition of $k$ before \er{f1}).
Then the identities \er{T34-1}, \er{T34-2} imply \er{T34-3}.
$\BBox$

\section {Proof of the Main Theorems}
\setcounter{equation}{0}

\noindent
{\bf Proof of Theorem \ref{T1}}. i) The asymptotics \er{T1-1}  were proved in Lemma \ref{T21}.

ii) and iii) of Lemma \ref{T32} give \er{st2} for the case
of non-virtual states, i.e., those not coinciding with $e_n^\pm$.

Lemma \ref{T33} implies \er{st2} for the case
of virtual states.

Lemma \ref{T311} gives \er{st1} for the case
of non-virtual states.
Lemma \ref{T33} implies  \er{st1} for the case
of virtual states.

ii) Using ii) and iii) of Lemma \ref{T32} we obtain \er{ajf}.

iii) Due to i) $\z$ is a zero of $F$, so \er{T22-3} yields \er{T1-2}.
 Lemma \ref{T22} ii) completes the proof of iii).
\BBox

\medskip
\no {\bf Proof of Theorem \ref{T2}}.
i) Let $g_n\ne \es$. The entire function $F=\vp(1,\cdot)\P_0^+\P_0^-$
has opposite sign on $\s_n$ and $\s_{n+1}$, since
$\P_0^+(z)\P_0^-(z)=|\P_0^+(z)|^2>0$ for $z$ inside $\s_n\cup \s_{n+1}$
(see \er{T21-1}) and $\vp(1,\cdot)$ has one simple zero in each interval
$[e_n^-,e_n^+]$.
Therefore $F$ has an odd number of zeros on $[e_n^-,e_n^+]$.

By Lemma \ref{T311}-\ref{T33}, $\z\in g_n^c$ is a state of $H$ if and only if
$\z\in \ol g_n$ is a zero of $F$ (with equal multiplicities).
Hence the number of states (counting multiplicities) on $g_n^c$ is odd.

Using Lemma \ref{T22} and \ref{T311}, we deduce that there exists exactly one
simple state $z_n$ in each interval $[e_n^-,e_n^+]$  for $g_n\ne \es$ and for
sufficiently large $n\ge 1$. Moreover, the asymptotics \er{ape}
$e_n^\pm=\pi n+{p_0\/2\pi n}+o(1/n)$ give
\[
\lb{rae}
z_n=\pi n+{p_0\/2\pi n}+o(1/n).
\]
By arguments analogous to the proof of \er{T21-2}, we obtain the identities
\[
\lb{T21-20}
\wt \vt(0,z)=1+\int_0^t\vp(x,z)q(x)\wt \vt(x,z)dx,\qqq
\wt \vp(0,z)=\int_0^t\vp(x,z)q(x)\wt \vp(x,z)dx;
\]
then a standard iteration procedure and \er{rae} give the asymptotics
\[
\lb{asvt1}
\wt \vt(0,z_n)=1+O(1/n),
\]
\[
\lb{asvp1}
\wt\vp(0,z_n)=\int_0^t{\sin^2 z_nx\/z_n^2}q(x)\,dx+O(1/n^3)
={q_0-\wh q_{cn}\/2(\pi n)^2}+O({1/n^3}),
\]
where $\wh q_{cn}=\int_0^tq(x)\cos 2\pi nx\,dx$.
Using \er{LD0} and $\P^\pm=\wt\vt+m_\pm\wt\vp$, see \er{Ptp},
 we obtain
\[
\lb{idF1}
F=F_1+F_2+F_3,\qqq F_1=\vp(1,\cdot)\wt\vt_0^2,\qq F_2=2\b\wt\vt_0\wt\vp_0,\qqq F_3=-\vt'(1,\cdot)\wt\vp_0^2,
\]
where we abbreviate $\wt\vt_0=\wt\vt(0,z),\ \wt\vp_0=\wt\vp(0,z)$.
Using the estimates \er{asvt1}, \er{asvp1}, we obtain
\begin{multline}
\lb{T31-4}
\qqq\qqq
F_1(z_n)=\vp(1,z_n)(1+O(1/n)),\qqq \qqq F_3(z_n)=O(1/n^{4})
\\
F_2(z_n)=(-1)^n{(p_{sn}+O({1/n}))\/\pi n} {(q_0-\wh q_{cn}+O({1/n}))\/2(\pi n)^2}=f_n+O(1/n^{4})\qq \as \ n\to \iy,
\qqq\qqq
\end{multline}
where $f_n=(-1)^n{p_{sn}(q_0-\wh q_{cn})\/2(\pi n)^3}$.
Combining these asymptotics with $F(z_n)=0$, we get
\[
\vp(1,z_n)=-F_2(z_n)+O(1/n^{4})=-f_n+O(1/n^{4}).
\]
Then, using $\vp(1,z_n)=\pa_z\vp(1,\m_n)\d_n+O(1/n^{4})$,
where $z_n=\m_n+\d_n$, we obtain
$$
\pa_z\vp(1,\m_n)\d_n=-f_n+O(1/n^{4})
$$
and the asymptotics $\pa_z\vp(1,\m_n)={(-1)^n\/(\pi n)}+O({1/n^2})$ give
$$
\d_n=-{f_n\/\pa_z\vp(1,\m_n)}+O(1/n^{3})
=-{(\wh q_0-\wh q_{cn})p_{sn}\/2(\pi n)^2}+O({1/n^3}),
$$
which yields \er{T2-1}.

Now we denote by $\cN (r,f)$ the total number of zeros of $f$ of modulus
$\le r$, and by $\cN^+(r,f)$ the number of zeros counted in $\cN (r,f)$
with non-negative real part, by $\cN^-(r,f)$ those with negative real part
(each zero being counted according to its multiplicity).
Then the following result holds [Koo].

\no    {\bf Theorem (Levinson).}
{\it  Let the entire function $f\in Cart_\o$.
Then $  \cN^{\pm }(r,f)={r\/ \pi }(\o+o(1))$ as $r\to \iy$,
and for each $\d >0$ the number of zeros of $f$ of modulus $\le r$
lying outside both of the two sectors $|\arg z | <\d$, $|\arg z -\pi |<\d$
is $o(r)$ for $r\to \iy$.}

We also denote by $\cN_+(r,f)$ (or $\cN_-(r,f)$) the number of zeros of $f$
counted in $\cN (r,f)$ with non-negative (or negative) imaginary part,
each zero being counted according to its multiplicity.

Let $s_0=0$ and $\pm s_n>0, n\in \N$, be all real zeros of $F$
and let $n_0$ be the multiplicity of the zero $s_0=0$.
Consider the entire function
$F_1=z^{n_0}\lim\limits_{r\to \iy}\prod\limits_{0<s_n\le r}(1-{z^2\/ s_n^2})$.
The Levinson Theorem and Lemma \ref{T21} imply
\[
\cN(r,F)=\cN(r,F_1)+\cN(r,F/F_1)=2r{1+2t\/\pi}+o(r),\qq
\cN(r,F_1)={2r\/\pi}+o(r),
\]
\[
\cN_-(r,F)=\cN_+(r,F)=\cN_-(r,\P_0^+)-N_0,
\]
as $r\to\iy$, where $N_0$ is the number of non-positive eigenvalues of $H$. Thus
\[
2\cN_-(r,F)=2r{2t\/\pi}+o(r),
\]
which yields \er{T2-2}.

ii) Using Lemma \ref{T34} we obtain the statements ii) and iii).
\BBox

\medskip
\noindent
{\bf Proof of Theorem \ref{T3}.}
Let $z=e_n^\pm$. Identity \er{T31-1} and $k(e_n^\pm)=\pi n$ yield
\[
\lb{Jep}
\P_0^-(z)=\P_0^+(z)=(-1)^Nw_+(z),\ \ \  w_+(z)=\F'(n_t,z)-{\b(z)\/\vp(1,z)}\F(n_t,z),\qq
N=n_tn.
\]
 Estimates \er{efs1} and $e_n^\pm=\pi n+\ve_n(p_0\pm |p_n|+O(\ve_n)), \ \ve_n={1\/2\pi n}$ (see \er{ape}), give
\[
\F'(n_t,z)=(-1)^N+{O(1/n)},\qq  \F(n_t,z)={\sin n_tz\/\pi n}+{O(1/n^2)}
={O(1/n^2)}.
\]
Using \er{efs}, we obtain
$$
\sin e_n^\pm=(-1)^n\sin { \pm |p_n|\/2\pi n}+O({1/n^2})
=(-1)^n{ \pm|p_n|\/2\pi n}+O({1/n^2}),
$$$$
\vp(1,e_n^\pm)={\sin e_n^\pm\/\pi n}+{(-1)^np_{cn}\/2\pi^2n^2} +{O(1/n^3)}=
(-1)^n{\pm|p_n|+p_{cn}\/2\pi^2n^2}+O({1/n^3}).
$$
Then the estimate $\sqrt{x^2+y^2}-y\ge x$ for $y,x\ge 0$ gives
$|p_n|\pm p_{cn}\ge |p_{sn}|$, which yields
\[
\lb{bdvp}
{\b(e_n^+)\/\vp(1,e_n^+)}=\pi n{p_{sn}+O({1/n})\/\pm|p_n|+p_{cn}+O({1/n})}=
O(\pi n),
\]
since $|p_{sn}|\ge {1\/n^\a}$. Combining \er{Jep}-\er{bdvp} and \er{abn}, we obtain
$\P_0^+(e_n^\pm)=1+o(1)$. The function $\P_0^+(z)$ is analytic on $g_n^-$
and $\P_0^+(e_n^\pm)=1+o(1)$. Thus $\P_0^+(z)$ has no zeros on $g_n^-$,
since by Theorem \ref{T2}, the function $F$ has exactly one zero on each
non-empty $\ol g_n$ for large $n>1$.

Let $\m_n+i0\in g_n^+$ be a bound state of $H_0$ for some sufficiently large
$n$.
Then Lemma \ref{Tm} implies $h_{sn}>0$. Moreover,
\er{anc} gives $h_{sn}=-{p_{sn}+O({1\/n})\/2\pi n}$ as $n\to \iy$.
Thus $p_{sn}<-{1\/n^\a}$ for large $n>1$ and, by the asymptotics \er{T2-1},
the bound state  $z_n>\m_n$ if $q_0>0$ and $z_n<\m_n$ if $q_0<0$.
The proof of the other cases is similar.
\BBox

\medskip
\noindent
{\bf Proof of Theorem \ref{T4}.} i)
 Using the identities \er{T21-2} and \er{Ptp} we obtain
\[
\lb{iii}
\P_0^+=Y_1+{i\sin k\/\vp_1}\wt \vp(z_n),\qq Y_1=\wt \vt(z_n)+{\b\/\vp_1}\wt \vp(z_n).
\]
Note that \er{LD0} gives $\b(\m_n)=0$. Then the asymptotics \er{T2-1}, \er{asb},
\er{efs1} imply
\[
{\b(z_n)\/\vp(1,z_n)}={\b'(\m_n)+O({\ve_n^3})\/\pa_z\vp(1,\m_n)+O({\ve_n^3})}
=o(1)\qqq \as \qq n\to \iy,\qqq \ve_n={1\/2\pi n}
\]
where we used the asymptotics $\pa_z\vp(1,\m_n)={(-1)^n+O({\ve_n})\/(\pi n)}$
and $\b'(\m_n)={o(1/n)}$. Thus \er{asvt1}, \er{asvp1} give
\[
Y_1=1+O(\ve_n), \qqq \wt\vp(0,z_n)=2\ve_n^2(b_n+O(\ve_n)), \qqq b_n=q_0-\wh q_{cn} 
\]
In the following we need these identities and asymptotics as $n\to \iy$
from \cite{KK}:
\[
\lb{35} (-1)^{n+1}i\sin k(z)=\sinh v(z)=\pm |\D^2(z)-1|^{1\/2}>0\qq
\ (z\in  g_n^\pm),
\]
\[
\lb{pav}
v(z)=\pm |(z-e_n^-)(e_n^+-z)|^{1\/2}(1+O(1/n^{2})),\qq
\sinh v(z)=v(z)(1+O(|g_n|^2),\qq (z\in \ol g_n^\pm).
\]
We rewrite the equation $\P_0^+=0$ in the form $\vp_1Y_1=-i\sin k\wt \vp(z_n)$.
Then we obtain
$$
2\d \ve_n(1+O(\ve_n))=v(z)2\ve_n^2(b_n+O(\ve_n))=\sqrt{\d(|g_n|-\d)}2\ve_n^2(b_n+O(\ve_n)),
$$
$$
\sqrt \d =\ve_n\sqrt{|g_n|-\d}(b_n+O(\ve_n)), \qqq \d=z_n-\m_n,
$$
where $\sqrt \d>0$ if $b_n>0$  and $\sqrt \d<0$ if $b_n<0$. The last
asymptotics imply $\d =\ve_n^2|g_n|(b_n+O(\ve_n))^2$, where
$b_n=q_0-\wh q_{cn}$, which yields \er{T4-1}.
\BBox

{\bf Proof of Theorem \ref{TIP}.}
i) Let $q\in \cQ_t, q_0=0$ and let each  $|\wh q_{cn}|>n^{-\a}$
for some $\a\in (0,1)$ and for sufficiently large $n$.
(The proof of the other cases is similar.)
Using the inverse spectral theory from \cite{K5} summarised in the Introduction above, we see that for any sequence  $\vk=(\vk_n)_1^\iy\in \ell^2, \vk_n\ge 0$,
there exists a potential $p\in L^2(0,1)$ such that gap $\g_n$ has length
$|\g_n|=\vk_n$ $(n\in\N)$. 
Moreover, for large $n$ we can write the gap in the form $\g_n=(E_n^-,E_n^+)$,
where $\m_n^2=E_n^-$ or $\m_n^2=E_n^+$.
We choose $\m_n^2=E_n^-$ or $\m_n^2=E_n^+$ depending on the given sequence
$\s=(\s_n)_{1}^\iy$, where $\s_n\in \{0,1\}$, as follows, using Theorem
\ref{T3} (i).

If $\s_n=1$ and $\wh q_{cn}<-n^{-\a}$ (or $\wh q_{cn}>n^{-\a}$), then taking
$\m_n^2=E_n^-$ (or $\m_n^2=E_n^+$) we ensure that $\l_n$ is an eigenvalue
for large $n$.

If $\s_n=0$ and $\wh q_{cn}>n^{-\a}$ (or $\wh q_{cn}<-n^{-\a}$), then taking
$\m_n^2=E_n^-$ (or $\m_n^2=E_n^+$) we ensure that $\l_n$ is an antibound state
for large $n$.

ii) Let $p\in L^2(0,1)$ such that $\m_n^2=E_n^-$ or $\m_n^2=E_n^+$ for
sufficiently large $n\in \N_0$, $\g_n=(E_n^-,E_n^+)$.
Let $\s=(\s_n)_{1}^\iy$ be any sequence, where $\s_n\in \{0,1\}$.
We take $|\wh q_{cn}|>n^{-\a}$ for large $n\in \N_0$, and then choose
the sign of $q_{cn}$ as follows, using Theorem \ref{T3} i).

If $\s_n=0$ and $\l_n^0=E_n^-$ (or $\l_n^0=E_n^+$), then taking
$\wh q_{cn}>n^{-\a}$ (or $\wh q_{cn}<-n^{-\a}$) we ensure that
$\l_n$ is an antibound state.

If $\s_n=1$ and $\l_n^0=E_n^-$ (or $\l_n^0=E_n^-$), then taking
$\wh q_{cn}<-n^{-\a}$ (or $\wh q_{cn}>n^{-\a}$) we ensure that
$\l_n$ is an eigenvalue.
\BBox

\medskip
\noindent
The proof of Theorem \ref{T5} will be based on the following asymptotic result
on the number of eigenvalues introduced in a gap by a compactly supported
perturbation. This follows the original idea of \cite{So} in the modified
version of \cite{Sc}.

\begin{lemma}
\lb{Sob}
Let $q \in \cQ_t$ be continuous and $H_\t$, $\O=[E_1, E_2]$ as in Corollary
\ref{T5}. Then
\[
\lb{asbs}
\#_{bs}(H_\t,\O)={\t}\int_0^\iy \rt(\r(E_2-q(x))-\r(E_1-q(x))\rt)dx+o(\t)\qq \as \qq \t\to \iy,
\]
where $\r$ is the integrated density of states \er{ids}.
\end{lemma}

\no {\bf Proof.}
By the Glazman decomposition principle and the properties of the periodic
problem on a half-line summarised in the Introduction above, the number of
bound states of $H_\t$ in $\O$ differs by no more than 5 from the number of
eigenvalues in $\O$ of the regular Sturm-Liouville operator $\tilde H_\t$ on
$[0, t\t]$ with Dirichlet boundary conditions. This number can be estimated by
oscillation theory as follows. For $\l\in\O$, let $u$ be the (real-valued)
solution of the equation
$$
  -u'' + p(x) u + q(x/\t) u = \l u
$$
with $u(0) = 0$, $u'(0) = 1$, and $\theta(\cdot, \l) : [0, t\t] \rightarrow
\R$ the locally absolutely continuous function (called a {\it Pr\"ufer
angle}) with the properties
$$
  \begin{pmatrix}u \cr u'\cr\end{pmatrix}=\sqrt{u^2 + u'^2}\,\begin{pmatrix}\sin\theta\cr\cos\theta\cr\end{pmatrix},
\qquad
  \theta(0,\l) = 0.
$$
Then $\theta(t\t,\l)$ is continuous and monotone increasing as a function of
$\l$, and the number of eigenvalues of $\tilde H_\t$ in $\O$ differs by no more
than 1 from $\frac{1}{\pi} (\theta(t\t,E_2)-\theta(t\t,E_1))$.

In the following we shall use the fact, which can be shown along the lines of
\cite{E} Theorem 3.1.2, \cite{Sc} Corollary 1, that the Pr\"ufer angle
$\tilde\theta$ of any real-valued solution of the periodic equation
$$
  -y'' + p y = \l y
$$
satisfies $\tilde\theta(x) = \pi \r(\l) x + O(1)$ as $x\rightarrow\infty$.

Consider a division of the interval $[0,t]$ into $n$ parts,
$0 = s_0 < s_1 < \dots < s_{n-1} < s_n = t$.
We set $q_j^+ := \sup\limits_{[s_{j-1}, s_j]} q$,
$q_j^- := \inf\limits_{[s_{j-1}, s_j]} q$.
Then if $\theta_j^\pm$ is the Pr\"ufer angle of a non-trivial real-valued
solution of
$$
  -y'' + p(x) y + q_j^\pm y = \l y
$$
on the interval $[\t s_{j-1}, \t s_j]$ with
$\theta_j^\pm(\t s_{j-1}) = \theta(\t s_{j-1})$, then Sturm comparison
(\cite{W} Thm 13.1) shows
that
$$
\theta_j^+(\t s_j) \le \theta(\t s_j) \le \theta_j^-(\t s_j).
$$
Hence
\begin{multline}
\lb{stuco}
\sum_{j=1}^n \r(\l-q_j^+)(s_j-s_{j-1})
= \lim_{\t\rightarrow\infty} \sum_{j=1}^n
  \frac{\theta_j^+(\t s_j,\l) - \theta_j^+(\t s_{j-1},\l)}{\t\pi}
\le \lim_{\t\rightarrow\infty} \sum_{j=1}^n
  \frac{\theta(\t s_j,\l) - \theta(\t s_{j-1},\l)}{\t\pi} \\
\le \lim_{\t\rightarrow\infty} \sum_{j=1}^n
  \frac{\theta_j^-(\t s_j,\l) - \theta_j^-(\t s_{j-1},\l)}{\t\pi}
= \sum_{j=1}^n \r(\l-q_j^-)(s_j-s_{j-1}).
\end{multline}
The extremes in \er{stuco} are Riemann sums of the integral
$\int\limits_0^t \r(\l-q(s))\,ds$; thus refining the division and observing
that $\r(\l-q(s)) = \r(\l) =$ const for all $\l\in\O$ if $s > t$, we obtain
$$
  \lim_{\t\rightarrow\infty} \frac{\theta(t\t,E_2)-\theta(t\t,E_1)}{\t\pi}
= \int_0^\infty (\r(E_2-q(s)) - \r(E_1-q(s)))\,ds,
$$
and the assertion follows.
\BBox

\medskip
\noindent
To complete the proof of Corollary \ref{T5}, we observe that
Theorem \ref{T2} (iii) implies that
$\#_{abs}(H_\t,\O^{(2)})\ge 1+\#_{bs}(H_\t,\O)$,
which together with \er{asbs} yields \er{T5-1}.

\bigskip
\no {\bf Acknowledgments.}
\setlength{\itemsep}{-\parskip} \footnotesize
Various parts of this paper were written at Universit\'e de Gen$\grave{\rm e}$ve, Section de Mathematiques and
Mathematical Institute of the Tsukuba Univ., Japan.
E. K. is grateful to these Institutes for their hospitality;
he would also like to thank A. Sobolev (London) for useful discussions about
the asymptotics associated with Theorem \ref{T5}. He was supported by the
14.740.11.0581 Federal Program "Development of Scientific Potential of Higher
Education for the period 2009--2013".
Both authors thank the referee for helpful detailed comments.



\begin{thebibliography}
{99999}\setlength{\itemsep}{-\parskip} \footnotesize

\bibitem[BZ] {BZ} Bindel, D.; Zworski, M.  Symmetry of Bound and Antibound States in the Semiclassical Limit,
Letters in Mathematical Physics,
81( 2007), 107-117.


\bibitem[BKW] {BKW}  Brown, B.; Knowles, I.; Weikard, R. On the
inverse resonance problem, J. London Math. Soc. (2) 68 (2003), no. 2, 383--401.


\bibitem[BW] {BW}
Brown, B.; Weikard, R. The inverse resonance problem for perturbations of algebro-geometric potentials. Inverse Problems 20 (2004), no. 2, 481--494.


\bibitem[Ch] {Ch} Christiansen, T. Resonances for steplike potentials: forward and inverse results. Trans. Amer. Math. Soc. 358 (2006), no. 5, 2071--2089.

\bibitem[CL] {CL}    Coddington, E.; Levinson, N. Theory of ordinary differential equations. New York, Toronto, London: McGraw-Hill 1955.


\bibitem[D] {D} Dimassi, M. Spectral shift function and
resonances for slowly varying perturbations of periodic
Schr\"odinger operators. J. Funct. Anal. 225 (2005), no. 1,
193--228.


\bibitem[DG] {DG} Dyatlov, S.; Ghosh, S. Symmetry of bound and antibound states in the semiclassical limit for a general class of potentials,
Proceeding of AMS,
 138(2010), No 9, 3203--3210.

\bibitem[E1] {Ea} Eastham, M.S.P. Theory of ordinary differential equations.
Van Nostrand Reinhold, London, 1970.

\bibitem[E2] {E} Eastham, M.S.P. The spectral theory of periodic differential equations. Scottish Academic Press, Edinburgh, 1973.

\bibitem[F1] {F1}  Firsova, N. Resonances of the perturbed Hill operator with exponentially decreasing extrinsic potential. Mat. Zametki 36 (1984), 711--724.

\bibitem[F2] {F2}     Firsova, N. The Levinson formula for a perturbed Hill operator. (Russian) Teoret. Mat. Fiz. 62 (1985), no. 2, 196--209.

\bibitem[F3] {F3}    Firsova, N. A direct and inverse scattering problem for a one-dimensional perturbed Hill operator.  Mat. Sb. 130(172) (1986), no. 3, 349--385.

\bibitem[F4] {F4}     Firsova, N. A trace formula for a perturbed one-dimensional Schr\"odinger operator with a periodic potential. I. (Russian) Problems in mathematical physics, No. 7 (Russian), pp. 162--177, Izdat. Leningrad. Univ., Leningrad, 1974.


\bibitem[Fr] {Fr}
Froese, R. Asymptotic distribution of resonances in one
dimension, J. Diff. Eq., 137( 1997), 251-272.

 \bibitem[GT]{GT} Garnett, J.; Trubowitz, E. Gaps and bands of
one dimensional  periodic Schr\"odinger operators II. Comment. Math.
Helv. 62(1987), 18-37.


\bibitem[G] {G} G\'erard, C. Resonance theory for periodic Schr\"odinger operators,
Bulletin de la S.M.F., 118(1990), n0 1, 27-54.

\bibitem[GS] {GS} Gesztesy, F.; Simon, B. A short proof of Zheludev's theorem. Trans. Amer. Math. Soc. 335 (1993), no. 1, 329--340.

\bibitem[HKS] {HKS} Hinton,D. ; Klaus, M.; Shaw, J. On the
Titchmarsh-Weyl function for the half-line perturbed periodic Hill's
equation. Quart. J. Math. Oxford Ser. (2) 41 (1990), no. 162,
189--224.

\bibitem[H] {H} Hitrik, M. Bounds on scattering poles in one dimension. Comm. Math. Phys. 208 (1999), no. 2, 381--411.

\bibitem[Ho] {Ho} Hochstadt, H. On the determination of a Hill's equation from its spectrum. Arc. Rat. Mech. Anal. 19 (1965) 353--362

\bibitem[IM] {IM} Its, A.; Matveev, V. Schr\"odinger operators with the
finite-gap  spectrum and the N-soliton solutions of the Korteweg--de Vries
 equation. Teoret. Math. Phys. 23 (1975), 51--68.

\bibitem[J] {J} Javrjan, V. Some perturbations of self-adjoint operators. (Russian) Akad. Nauk Armjan. SSR Dokl. 38 (1964), 3--7.


\bibitem[KK] {KK}  Kargaev, P.; Korotyaev, E. Effective masses and
conformal mappings. Comm. Math. Phys. 169 (1995), no. 3, 597--625.

\bibitem[KK1]{KK1}
 Kargaev, P.; Korotyaev, E. Inverse Problem for the Hill
Operator, the Direct Approach.  Invent. Math., 129 (1997), no. 3, 567--593.


\bibitem[KM] {KM} Klopp, F.; Marx, M.
The width of resonances for slowly varying perturbations of one-dimensional periodic Schr\"odingers operators,  Seminaire: Equations aux Deriv\'ees Partielles. 2005--2006, Exp. No. IV, 18 pp., Semin. Equ. Deriv. Partielles, Ecole Polytech., Palaiseau.



\bibitem[K1] {K1} Korotyaev, E. Inverse resonance scattering  on the half
line. Asymptotic Anal.  37 (2004), No 3/4, 215--226.

\bibitem[K2] {K2} Korotyaev, E. Inverse resonance scattering on the real line. Inverse Problems 21 (2005), no. 1, 325--341.

\bibitem[K3] {K3} Korotyaev, E. Stability for inverse resonance problem.
Int. Math. Res. Not. 2004, no. 73, 3927--3936.

\bibitem[K4] {K4} Korotyaev, E. Resonance theory for perturbed Hill operator.
Asymptotic Anal. (to appear)

\bibitem[K5] {K5} Korotyaev, E. Inverse problem and the trace formula for the Hill operator. II Math. Z. 231 (1999), no. 2, 345--368.

\bibitem[K6] {K6}
  Korotyaev, E. Characterization of the spectrum of
 Schr\"odinger operators with periodic distributions.
 Int. Math. Res. Not. 2003, no. 37, 2019--2031.

 \bibitem[K7] {K7}
 Korotyaev, E. Estimates of periodic potentials in terms of gap
lengths. Comm. Math. Phys. 197 (1998), no. 3, 521--526.



\bibitem[Koo] {Koo} Koosis, P. The logarithmic integral I, Cambridge Univ.
Press, Cambridge, London, New York 1988.

\bibitem[M] {M} Marchenko, V. Sturm-Liouville operator and applications.
Basel: Birkh\"auser 1986.

\bibitem[MO] {MO}  Marchenko, V.; Ostrovski I. A characterization of
the spectrum of the Hill operator. Math. USSR  Sbornik  26(1975), 493-554.

\bibitem[N-Z] {N-Z} Novikov, S.; Manakov, S. V.; Pitaevskii, L. P.;
Zakharov, V. E. Theory of solitons. The inverse scattering method.
Translated from the Russian. Contemporary Soviet Mathematics.
Consultants Bureau [Plenum], New York, 1984.


\bibitem[PT] {PT} P\"oschel, P.; Trubowitz, E. Inverse Spectral Theory.
Boston: Academic Press, 1987.

\bibitem[Rb] {Rb}    Rofe-Beketov, F. A finiteness test for the number
of discrete levels which can be introduced into the gaps of the continuous spectrum by perturbations of a periodic potential. Dokl. Akad. Nauk SSSR 156 (1964), 515--518.

\bibitem[Sc] {Sc}  Schmidt, K. M.
Eigenvalue asymptotics of perturbed periodic Dirac systems in the slow-decay limit. Proc. Amer. Math. Soc. 131 (2003) 1205-1214.

\bibitem[S1] {Si} Simon, B. On the genericity of nonvanishing instability intervals in Hill's equation. Ann. Inst. Henri Poincar\'e 24 (1976), 91--93

\bibitem[S2] {S}  Simon, B. Resonances in one dimension and Fredholm
determinants. J. Funct. Anal. 178 (2000), no. 2, 396--420.

\bibitem[So] {So}  Sobolev, A.V.  Weyl asymptotics for the discrete
spectrum of the perturbed Hill operator. Estimates and asymptotics for discrete spectra of integral and differential equations (Leningrad, 1989--90), 159--178, Adv. Soviet Math., 7, Amer. Math. Soc., Providence, RI, 1991.



\bibitem[T] {T} Titchmarsh, E. Eigenfunction expansions associated
with second-order differential equations 2, Clarendon Press, Oxford, 1958.




\bibitem[Tr] {Tr} Trubowitz, E. The inverse problem for  periodic potentials.
Commun. Pure Appl. Math. 30(1977), 321-337.

\bibitem[W] {W} Weidmann, J. Spectral theory of ordinary differential operators. Lect. Notes in Math. 1258, Springer, Berlin, New York 1987


\bibitem[Z] {Z}  Zworski, M. Distribution of poles for scattering on
the real line, J. Funct. Anal. 73(1987), 277-296.

\bibitem[Z1] {Z1} Zworski, M. SIAM, J.  Math. Analysis, "A remark on
isopolar potentials" 82(6), 2002, 1823-1826.


\bibitem[Z2] {Z2} Zworski, M. Counting scattering poles. In: Spectral
and scattering theory (Sanda, 1992),
301--331, Lecture Notes in Pure and Appl. Math., 161, Dekker, New York, 1994.


\bibitem[Zh1] {Zh1}  Zheludev, V. A. The eigenvalues of a
perturbed Schr\"odinger operator with periodic potential. (Russian) 1967 Problems of Mathematical Physics, No. 2, Spectral Theory, Diffraction Problems pp. 108--123.

\bibitem[Zh2] {Zh2}  Zheludev, V. A. The perturbation of the spectrum of the Schr\"odinger operator with a complex-valued periodic potential. (Russian) Problems of mathematical physics,  Spectral theory 3(1968),  31--48.

\bibitem[Zh3] {Zh3}  Zheludev, V. The spectrum of Schr\"odinger's operator,
with a periodic potential, defined on the half-axis. Works of Dept. of Math. Analysis of Kaliningrad State University (1969) (Russian), pp. 18--37.





\end{thebibliography}
\end{document}